\documentclass[reqno]{amsart}
\usepackage{amsfonts,latexsym,amssymb,amsmath,amsthm,color,enumerate}

\usepackage[pdftex]{graphicx}
\usepackage[all]{xy}
\usepackage{tikz}
\usetikzlibrary{shapes}
\usetikzlibrary{plotmarks}

\newtheorem{teorema}{Theorem}[section]
\newtheorem{propo}[teorema]{Proposition}
\newtheorem{lemma}[teorema]{Lemma}
\newtheorem{corollary}[teorema]{Corollary}

\newtheorem{Example}{Example}
\theoremstyle{definition}

\newtheorem{properties}[teorema]{Properties}
\newtheorem{property}[teorema]{Property}
\newtheorem{Remark}[teorema]{Remark}

\theoremstyle{remark}
\newtheorem{nota}[teorema]{Remark}

\newcommand{\ord}{{\rm{ord}}}
\newcommand{\Zer}{{\rm{Zer}}}
\newcommand{\Char}{{\rm{Char}}}
\newcommand{\supp}{{\hbox{\rm supp}}}
\newcommand{\cont}{{\hbox{\rm cont}}}
\newcommand{\ini}{{\hbox{\rm in}}}

\newcommand{\C}{\mathbb{C}}
\newcommand{\N}{\mathbb{N}}
\newcommand{\R}{\mathbb{R}}

\newcommand{\cN }{\mathcal{N}}

\newcommand{\Teisssr}[4]{\Bigl\{
   \setlength{\unitlength}{1ex}
   \begin{picture}(#3,3)(0,0.4)
      \put(0,1.15){\line(1,0){#3}}
      \put(0,0.85){\line(1,0){#3}}
      \put(#4,1.3){\makebox(0,0)[b]{$#1$}}
      \put(#4,0.7){\makebox(0,0)[t]{$#2$}}
   \end{picture}\Bigr\}}

\begin{document}

\title[On the equisingularity class of the general higher order polars]{On the equisingularity class of the general higher order polars of plane branches}

\author{Evelia R. Garc\'{i}a Barroso}
\address[Evelia R. Garc\'{i}a Barroso]{Dpto. Matem\'{a}ticas, Estad\'{\i}stica e Investigaci\'on Operativa. IMAULL.\\
Universidad de La Laguna. Apartado de Correos 456. 38200 La Laguna, Tenerife, Spain.}
\email{ergarcia@ull.es}

\author {Janusz Gwo\'zdziewicz}
\address[Janusz Gwo\'zdziewicz]{Institute of Mathematics\\
University of the National Education Commission, Krakow\\
Podchor\c a{\accent95 z}ych 2\\
PL-30-084 Cracow, Poland}
\email{janusz.gwozdziewicz@up.krakow.pl}

\author{Mateusz Masternak}
\address[Mateusz Masternak]{Institute of Mathematics \\
Faculty of Mathematics and Natural Sciences\\
Jan Kochanowski University \\
ul. \'Swi\c{e}tokrzyska 15A\\
PL 25-406 Kielce, Poland.}
\email{mateusz.masternak@ujk.edu.pl}

\thanks{The first two authors were partially supported by
the Spanish grant PID2019-105896GB-I00 funded by
MCIN/AEI/10.13039/501100011033.}

\date{\today}

\begin{abstract}
In this paper we describe the factorization of the higher order polars of a generic branch 
in its equisingularity class. We generalize the results of  Casas-Alvero 
and Hefez-Hernandes-Hern\'andez  to higher order polars. 
\end{abstract}

\maketitle
\section{Introduction}

Let $f(x,y)\in \C[[x,y]]$ be an irreducible formal power series and  $C=\{f(x,y)=0\}$ be the  {\it branch} determined by $f(x,y)=0$. 
The {\it multiplicity} of $C$ is the order of~$f$. 
When this multiplicity is $n>1$ we say that $C$ is 
{\it singular}.  Otherwise $C$ is a {\it smooth} branch. In this paper we will consider singular branches. After a change of coordinates, if necessary, we may assume that $x=0$ is not tangent to $C$ at $0$. This is equivalent to $\ord f(0,y)=\ord f=n$. By Newton Theorem there is $\alpha(x)=\sum_{i\geq n}a_ix^{i/n}\in \mathbb C[[x^{1/n}]]\subset \mathbb C[[x]]^{*}$  such that $f(x,\alpha(x))=0$, where $\mathbb C[[x]]^{*}$ denotes the ring of Puiseux power series. The power series $\alpha(x)$ is called a {\it Newton-Puiseux root} of $f(x,y)$. It is well-known that the set of all Newton-Puiseux roots of $f(x,y)$ is $\Zer f:=\{\alpha_{\epsilon}(x)=\sum_{i\geq n}a_i\epsilon^i x^{i/n}\;:\;\epsilon \in \mathbb U_{n}\}$, where 
$\mathbb U_{n}$ is the multiplicative group of 
$n$th complex roots of unity. By Puiseux Theorem

\begin{equation}
    \label{NPf}
f(x,y)=u(x,y)\prod_{\epsilon\in \mathbb U_n}(y-\alpha_{\epsilon}(x)),
\end{equation}
where $u(x,y)$ is a unit in $\mathbb C[[x,y]]$.

The {\it index} of any $\alpha \in \mathbb C[[x]]^{*}$ is the smallest natural number $m$ such that
$\alpha \in \mathbb C[[x^{1/m}]]$.
To any $\alpha(x)=\sum_{i}a_ix^{i/n}\in \mathbb C[[x]]^{*}$ of positive order and index $n$ we associate with two finite sequences $(e_{i})_{i}$ and $(b_{i})_{i}$ of natural numbers as follows: $e_{0}=b_{0}=n$; if $e_{k}> 1$ \label{page1}
then $b_{k+1}:=\min\{i\;:\; a_i\neq 0;\; \gcd(e_k,i)<e_k\}$ and 
$e_{k+1}:=\gcd(e_k, b_{k+1})$. 
 The sequence $(e_{i})_{i}$ is strictly decreasing and for some $h \in \mathbb N$ we have $e_h=1$. The sequence  $(b_{0},b_{1},\ldots,b_{h})$ is called the {\it characteristic} of $\alpha$. 
By \cite[Lemma 6.8]{Hefez} we get


\begin{equation}\label{eq:1} \ord (\alpha_{\epsilon}(x)-\alpha(x))=\tfrac{b_{j}}{n}\;\; 
\hbox{\rm if and only if } \epsilon \in \mathbb U_{e_{j-1}} \backslash \mathbb U_{e_{j}}.
\end{equation}

Let $\lambda_l(x)$ be the sum of all terms of $\alpha(x)$ of degree strictly less than $\frac{b_l}{b_0}$. We denote by $f_l(y)$ the minimal polynomial of $\lambda_l(x)$ in the ring $\mathbb C[[x]][y]$. The
polynomial $f_l(y)$ does not depend on the choice of $\alpha(x)\in \Zer f$ and its degree is $\frac{n}{e_{l-1}}$.\\

Observe that the characteristic of  $ \alpha_{\epsilon}$ equals the characteristic of $\alpha$. The {\it characteristic} of an irreducible power series $f(x,y)\in \mathbb C[[x,y]]$ is the characteristic of any of its Newton-Puiseux roots. The set of {\em characteristic exponents} of $f$ is $\Char(f)=\left \{\frac{b_{i}}{n} \;:\; i\in \{1,\ldots,h\}\right  \}$. After \eqref{eq:1} the characteristic exponents of $f$ are the orders of differences of any two of its distinct Newton-Puiseux roots. \\
 
 Let $C=\{f(x,y)=0\}$ and $D=\{g(x,y)=0\}$ be two curves with $f,g\in \mathbb C[[x,y]]$. The {\it intersection multiplicity} of $C$ and $D$ is $i_0(C,D)=\dim \mathbb C[[x,y]]/(f,g)$ where $(\cdot, \cdot)$ denotes the ideal generated by two power series. Usually $i_0(C,D)$ is also denoted by $i_0(f,g)$.
 
 If $C$ and $D$ are branches then the {\em contact} of $C$ and $D$ is
 \[
 \cont (C,D)=\cont (f,g)=\max \{\ord (\alpha-\gamma)\;:\;\alpha\in \Zer \,f,\;\gamma \in \Zer \,g\}.
 \]

If $\alpha$ is a Puiseux series and $v\in \mathbb C[[x,y]]$ is irreducible then we put
\[
 \cont (\alpha,v)= \max \{\ord (\alpha-\gamma)\;:\;\gamma \in \Zer \,v\}.
 \]
 We say that the branches $C$ and $D$ are {\it equisingular} if and only if they have the same characteristic. We will denote by 
$K(b_{0},b_{1},\ldots,b_{h})$
 the coset of equisingular branches of characteristic  $(b_{0},b_{1},\ldots,b_{h})$. If $C=\{f(x,y)=0\}$ is a branch in $K(b_{0},b_{1},\ldots,b_{h})$,  by abuse of language we will put $f\in K(b_{0},b_{1},\ldots,b_{h})$. Let $f(x,y)=\sum_{ij}a_{ij}x^{i}y^{j}\in K(b_{0},b_{1},\ldots,b_{h})$. We can consider $f$ as an element of the  polynomial ring $\mathbb C[a_{ij}]$ with infinitely many  variables $a_{ij}$.
 We say that $f\in K(b_{0},b_{1},\ldots,b_{h})$ is {\it generic} in its equisingularity class if 
 its coefficients $a_{ij}$ belong to an open Zariski set.

\begin{figure} 
\begin{center}
\begin{tikzpicture}[x=0.7cm,y=0.7cm] 
\fill[fill=orange!20!white] (0,3) --(0,2)--(3,0)-- (4,0) --(4,3)--cycle;
\draw[->] (0,0) -- (4,0);
\draw[->] (0,0) -- (0,3);
\draw[-, line width=0.5mm] (0,2) -- (3,0);
\node [left] at (0,2) {$n$};
\node [below] at (3,0) {$m$};
\node [below] at (2.6,2) {$\Teisssr{m}{n}{2.2}{1.1}$};
\end{tikzpicture}
\end{center}
   \label{fig:elementary polygon}
   \caption{Elementary Newton diagram}
    \end{figure}
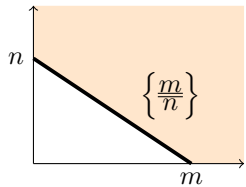

Let $A$ be a nonempty subset of $\N\times \N$. The {\it Newton diagram} $\cN(A)$ of the set $A$ is the convex hull of
$A+(\R_{\geq 0})^{2}$, where $+$ means the Minkowski sum. By definition, the {\em support} of any Newton diagram $\Delta$ is $\supp(\Delta):=\Delta\cap \N^2$.
We say that  $\cN(A)$ is {\it convenient} if it intersects both coordinate axes. The {\it Newton polygon} of the Newton diagram $\Delta$ is the union of the compact edges  of the boundary of $\Delta$, and we will denote it by $\delta(\Delta)$. A convenient Newton diagram is {\it elementary} if its boundary has exactly one compact edge. In this case, following  Teissier \cite{Teissier1}, we will denote by $\Teisssr{m}{n}{3}{1.5}$  the elementary Newton diagram of $A=\{(m,0), (0,n)\}$, for any positive natural numbers $m,n$ 
(see Figure \ref{fig:elementary polygon}).

 The {\it inclination} of the elementary Newton diagram $\Teisssr{m}{n}{3}{1.5}$ (and of any of its translations) is $m/n$.
Any convenient Newton diagram $\mathcal N$ can be written as a Minkowski sum of elementary Newton diagrams, where inclinations of successive elementary diagrams form an strictly decreasing sequence. This writing is called the {\it canonical representation} of  $\mathcal N$. A Newton diagram $\mathcal N$ can also be written as a sum of elementary Newton diagrams $\mathcal N=\sum_{i=1}^{r}\Teisssr{m_{i}}{n_{i}}{3}{1.5}$
where $\gcd(m_{i},n_{i})=1$ for any $i\in \{1,\ldots,r\}$ and $m_{i}/n_{i}\geq m_{i+1}/n_{i+1}$ for $i\in \{1,\ldots,r-1\}$. This new writing is called the {\it long canonical representation} of $\mathcal N$. The long canonical representation is unique.
\begin{Example}
The long canonical representation of $\mathcal N=\Teisssr{10}{4}{3}{1.5}+\Teisssr{8}{6}{3}{1.5}$ is \[\mathcal N=\Teisssr{5}{2}{3}{1.5}+\Teisssr{5}{2}{3}{1.5}+\Teisssr{4}{3}{3}{1.5}+\Teisssr{4}{3}{3}{1.5}.\]

In Figure \ref{xxx} we illustrate both canonical representations.\\ 
If we drop the hypothesis of $\gcd(m_{i},n_{i})=1$  in the definition of the long canonical representation  we can express $\mathcal N$ in other ways that are not canonical, for example
\[\mathcal N=\Teisssr{10}{4}{3}{1.5}+\Teisssr{4}{3}{3}{1.5}+\Teisssr{4}{3}{3}{1.5}.\]

\begin{figure}[h!] 
\begin{center}
\begin{tikzpicture}[x=0.35cm,y=0.35cm] 
\tikzstyle{every node}=[font=\small]
\fill[fill=orange!20!white] (0,11) --(0,10)--(8,4)-- (18,0) --(19,0)--(19,11)--cycle;
\foreach \x in {0,...,19}{
\foreach \y in {0,...,11}{
       \node[draw,circle,inner sep=1pt,fill, color=gray!40] at (1*\x,1*\y) {}; }
   }
\foreach \y in {1,2,...,10} \draw(0,\y)node[left]{\y};
\foreach \x in {0,1,...,18} \draw(\x,0)node[below]{\x};
 \foreach \y in {4,...,10,11}{
\foreach \x in {8,...,19}{
       \node[draw,circle,inner sep=1pt,fill, color=black] at (1*\x,1*\y) {}; }
   }
   
    \foreach \y in {2,3}{
\foreach \x in {13,...,19}{
       \node[draw,circle,inner sep=1pt,fill, color=black] at (1*\x,1*\y) {}; }
   }
      \foreach \y in {7,...,11}{
\foreach \x in {4,...,7}{
       \node[draw,circle,inner sep=1pt,fill, color=black] at (1*\x,1*\y) {}; }
   }
   \foreach \y in {7,...,11}{
\foreach \x in {4,...,7}{
       \node[draw,circle,inner sep=1pt,fill, color=black] at (1*\x,1*\y) {}; }
   }
    \foreach \y in {10,11}{
\foreach \x in {0,...,4}{
       \node[draw,circle,inner sep=1pt,fill, color=black] at (1*\x,1*\y) {}; }
   }
   \foreach \y in {1}{
\foreach \x in {16,...,19}{
       \node[draw,circle,inner sep=1pt,fill, color=black] at (1*\x,1*\y) {}; }
   }
\node[draw,circle,inner sep=1pt,fill, color=black] at (2,9){};
\node[draw,circle,inner sep=1pt,fill, color=black] at (3,9){};
\node[draw,circle,inner sep=1pt,fill, color=black] at (3,8){};
\node[draw,circle,inner sep=1pt,fill, color=black] at (6,6){};
\node[draw,circle,inner sep=1pt,fill, color=black] at (7,5){};
\node[draw,circle,inner sep=1pt,fill, color=black] at (7,6){};
\node[draw,circle,inner sep=1pt,fill, color=black] at (11,3){};
\node[draw,circle,inner sep=1pt,fill, color=black] at (12,3){};
\node[draw,circle,inner sep=1pt,fill, color=black] at (18,0){};
\node[draw,circle,inner sep=1pt,fill, color=black] at (19,0){};

\node [below] at (4,4){$8$};
\node [left] at (8,1.7){$4$};
\node [above] at (2,7){\textcolor{blue}{$4$}};
\node [right] at (4,5.5){\textcolor{blue}{$3$}};
\node [below] at (10.5,2){\textcolor{blue}{$5$}};
\node [left] at (13,1){\textcolor{blue}{$2$}};
\draw[->] (0,0) -- (19,0);
\draw[->] (0,0) -- (0,11);
\draw[-, line width=0.5mm] (0,10) -- (8,4);
\draw[-, line width=0.5mm] (18,0) -- (8,4);
\draw[-, dashed, line width=0.3mm] (0,4) -- (8,4);
\draw[-, dashed, line width=0.3mm] (8,0) -- (8,4);
\draw[-, dashed, line width=0.3mm, color=blue] (0,7) -- (4,7);
\draw[-, dashed, line width=0.3mm, color=blue] (4,4) -- (4,7);
\draw[-, dashed, line width=0.3mm, color=blue] (8,2) -- (13,2);
\draw[-, dashed, line width=0.3mm, color=blue] (13,0) -- (13,2);
\end{tikzpicture}
\end{center}
\caption{Canonical and long canonical representation of $\Teisssr{10}{4}{3}{1.5}+\Teisssr{8}{6}{3}{1.5}$}
\label{xxx}
    \end{figure}
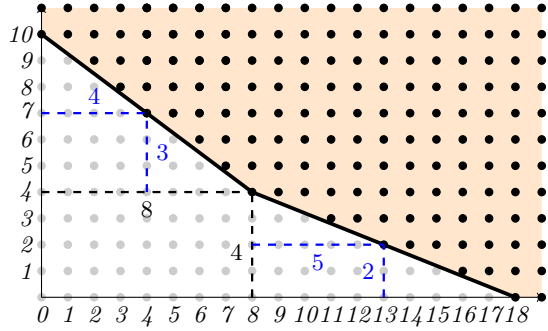
\end{Example}

The Newton diagram $\cN(f)$ of a nonzero power series $f(x,y)=\sum_{i,j}a_{ij}x^{i}y^{j}$ is the Newton diagram
$\cN(\supp(f))$, where  $\supp(f):=\{(i,j)\in \N^{2}\;:\; a_{ij}\neq 0\}$ is the support of $f$. 
It is well-known (see \cite[Lemme 8.4.2]{Chenciner}) that if $\mathcal \sum_{i=1}^{r}\Teisssr{M_{i}}{N_{i}}{3}{1.5}$ 
is the canonical representation of $\cN(f)$ then for any $i\in \{1,\ldots,r\}$ 
there are exactly $N_i$ Newton-Puiseux roots of $f$ of order $\frac{M_i}{N_i}$. 
Let $S$ be a compact edge of $\cN(f)$ of inclination $p/q$, where $p$ and $q$ are
coprime integers. The initial part of $f(x,y)$ with respect to $S$ is the 
quasi-homogeneous polynomial $f_S(x,y) = \sum_{ij} a_{ij}x^iy^j$ where the sum runs over all points in $S\cap \supp(f)$. Let $f_S(x,y)=ax^ky^l\prod_{j=1}^r(y^q-c_jx^p)^{s_j}$ be the factorization of $f_S$ into irreducible factors, where $k,l$ are non-negative integers and $a,c_j\in \mathbb C\backslash\{0\}$ with $c_j$ pairwise different.
The power series $f(x,y)$ is {\it non-degenerate} (in the sense of Kouchnirenko \cite{Kouchnirenko}) on  $S$ if one of the following equivalent conditions holds:
\begin{enumerate}
    \item[(ND1)] $s_j=1$ for any $j\in \{1,\ldots,r\}$.
    \item[(ND2)] All non-zero roots of $f_S(1,y)$ are simple.
    \item[(ND3)] All Newton-Puiseux roots of $f$ of order $p/q$ have different initial coefficients.
\end{enumerate} 


Let $\Delta$ be a Newton diagram and $k$ a nonnegative integer. 
The {\em symbolic $k$th derivative} $\Delta^{(k)}$ of $\Delta$ is the Newton diagram of the set $(\Delta-(0,k))\cap \N^2$.

\begin{figure}[h!] 
\begin{center}
\begin{tikzpicture}[scale=0.4]
\begin{scope}[shift={(10,10)}]
\tikzstyle{every node}=[font=\small]
\fill[fill=orange!20!white] (0,6) --(0,5)-- (12,0) --(13,0)--(13,6)--cycle;
\foreach \x in {0,...,12}{
\foreach \y in {0,...,5}{
       \node[draw,circle,inner sep=1pt,fill, color=gray!40] at (1*\x,1*\y) {}; }
   }
    \foreach \x in {5,...,10,11,12}{
\foreach \y in {3,4,5}{
       \node[draw,circle,inner sep=1pt,fill, color=black] at (1*\x,1*\y) {}; }
   }
    \foreach \x in {8,...,10,11,12}{
\foreach \y in {2}{
       \node[draw,circle,inner sep=1pt,fill, color=black] at (1*\x,1*\y) {}; }
   }
   \foreach \x in {0,...,4}{
\foreach \y in {5}{
       \node[draw,circle,inner sep=1pt,fill, color=black] at (1*\x,1*\y) {}; }
   }
    \foreach \x in {3,4}{
\foreach \y in {4}{
       \node[draw,circle,inner sep=1pt,fill, color=black] at (1*\x,1*\y) {}; }
   }

    \foreach \x in {11,12}{
\foreach \y in {1}{
       \node[draw,circle,inner sep=1pt,fill, color=black] at (1*\x,1*\y) {}; }
   }
   
   \foreach \y in {1,2,...,5} \draw(0,\y)node[left]{\y};
\foreach \x in {0,1,...,12} \draw(\x,0)node[below]{\x};
\node[draw,circle,inner sep=1pt,fill, color=black] at (12,0){};
\node[draw,circle,inner sep=1pt,fill, color=black] at (10,1){};
\node [right] at (5,-2){$\Delta=\Teisssr{12}{5}{3}{1.5}$};
\draw[->] (0,0) -- (13,0);
\draw[->] (0,0) -- (0,6);
\draw[-, line width=0.5mm] (0,5) -- (12,0);

   \end{scope}

   \begin{scope}[shift={(0,0)}]
   

\foreach \x in {0,...,12}{
\foreach \y in {1,...,5}{
       \node[draw,circle,inner sep=1pt,fill, color=gray!40] at (1*\x,1*\y) {}; }
   }
    \foreach \x in {5,...,10,11,12}{
\foreach \y in {3,4,5}{
       \node[draw,circle,inner sep=1pt,fill, color=black] at (1*\x,1*\y) {}; }
   }
    \foreach \x in {8,...,10,11,12}{
\foreach \y in {2}{
       \node[draw,circle,inner sep=1pt,fill, color=black] at (1*\x,1*\y) {}; }
   }
   \foreach \x in {0,...,4}{
\foreach \y in {5}{
       \node[draw,circle,inner sep=1pt,fill, color=black] at (1*\x,1*\y) {}; }
   }
    \foreach \x in {3,4}{
\foreach \y in {4}{
       \node[draw,circle,inner sep=1pt,fill, color=black] at (1*\x,1*\y) {}; }
   }

    \foreach \x in {11,12}{
\foreach \y in {1}{
       \node[draw,circle,inner sep=1pt,fill, color=black] at (1*\x,1*\y) {}; }
   }
   
   \foreach \y in {1,2,...,4} \draw(0,\y+1)node[left]{\y};
\foreach \x in {0,1,...,12} \draw(\x,-0.2)node[above]{\x};
\node [right] at (3,-1){$(\Delta-(0,1))\cap \mathbb N\times \mathbb N$};
\node[draw,circle,inner sep=1pt,fill, color=black] at (10,1){};
\draw[->] (0,1) -- (13,1);
\draw[->] (0,1) -- (0,6);

\end{scope}
 \begin{scope}[shift={(20,0)}]

\fill[fill=orange!20!white] (0,6) --(0,5)-- (10,1) --(13,1)--(13,6)--cycle;
\foreach \x in {0,...,12}{
\foreach \y in {1,...,5}{
       \node[draw,circle,inner sep=1pt,fill, color=gray!40] at (1*\x,1*\y) {}; }
   }
    \foreach \x in {5,...,10,11,12}{
\foreach \y in {3,4,5}{
       \node[draw,circle,inner sep=1pt,fill, color=black] at (1*\x,1*\y) {}; }
   }
    \foreach \x in {8,...,10,11,12}{
\foreach \y in {2}{
       \node[draw,circle,inner sep=1pt,fill, color=black] at (1*\x,1*\y) {}; }
   }
   \foreach \x in {0,...,4}{
\foreach \y in {5}{
       \node[draw,circle,inner sep=1pt,fill, color=black] at (1*\x,1*\y) {}; }
   }
    \foreach \x in {3,4}{
\foreach \y in {4}{
       \node[draw,circle,inner sep=1pt,fill, color=black] at (1*\x,1*\y) {}; }
   }

    \foreach \x in {11,12}{
\foreach \y in {1}{
       \node[draw,circle,inner sep=1pt,fill, color=black] at (1*\x,1*\y) {}; }
   }
   
   \foreach \y in {1,2,...,4} \draw(0,\y+1)node[left]{\y};
\foreach \x in {0,1,...,12} \draw(\x,-0.2)node[above]{\x};
\node [right] at (5,-1.2){$\Delta^{(1)}=\Teisssr{10}{4}{3}{1.5}$};
\node[draw,circle,inner sep=1pt,fill, color=black] at (10,1){};
\draw[->] (0,1) -- (13,1);
\draw[->] (0,1) -- (0,6);
\draw[-, line width=0.5mm] (0,5) -- (10,1);
\end{scope}

\begin{scope}[shift={(0,-10)}]
\foreach \x in {0,...,12}{
\foreach \y in {2,...,5}{
       \node[draw,circle,inner sep=1pt,fill, color=gray!40] at (1*\x,1*\y) {}; }
   }
    \foreach \x in {0,...,10,11,12}{
\foreach \y in {5}{
       \node[draw,circle,inner sep=1pt,fill, color=black] at (1*\x,1*\y) {}; }
   }
    \foreach \x in {3,...,9,10,11,12}{
\foreach \y in {4}{
       \node[draw,circle,inner sep=1pt,fill, color=black] at (1*\x,1*\y) {}; }
   }
   \foreach \x in {5,...,12}{
\foreach \y in {3}{
       \node[draw,circle,inner sep=1pt,fill, color=black] at (1*\x,1*\y) {}; }
   }
    \foreach \x in {8,...,12}{
\foreach \y in {2}{
       \node[draw,circle,inner sep=1pt,fill, color=black] at (1*\x,1*\y) {}; }
   }

   \foreach \y in {1,2,...,3} \draw(0,\y+2)node[left]{\y};
\foreach \x in {0,1,...,12} \draw(\x,0.8)node[above]{\x};
\node [right] at (3,-0.5){$(\Delta-(0,2))\cap \mathbb N\times \mathbb N$};
\draw[->] (0,2) -- (13,2);
\draw[->] (0,2) -- (0,6);
\end{scope}

\begin{scope}[shift={(20,-10)}]
\fill[fill=orange!20!white] (0,6) --(0,5)--(5,3)-- (8,2) --(13,2)--(13,6)--cycle;
\foreach \x in {0,...,12}{
\foreach \y in {2,...,5}{
       \node[draw,circle,inner sep=1pt,fill, color=gray!40] at (1*\x,1*\y) {}; }
   }
    \foreach \x in {0,...,10,11,12}{
\foreach \y in {5}{
       \node[draw,circle,inner sep=1pt,fill, color=black] at (1*\x,1*\y) {}; }
   }
    \foreach \x in {3,...,9,10,11,12}{
\foreach \y in {4}{
       \node[draw,circle,inner sep=1pt,fill, color=black] at (1*\x,1*\y) {}; }
   }
   \foreach \x in {5,...,12}{
\foreach \y in {3}{
       \node[draw,circle,inner sep=1pt,fill, color=black] at (1*\x,1*\y) {}; }
   }
    \foreach \x in {8,...,12}{
\foreach \y in {2}{
       \node[draw,circle,inner sep=1pt,fill, color=black] at (1*\x,1*\y) {}; }
   }

   \foreach \y in {1,2,...,3} \draw(0,\y+2)node[left]{\y};
\foreach \x in {0,1,...,12} \draw(\x,0.8)node[above]{\x};
\node [right] at (3,-0.5){$\Delta^{(2)}=\Teisssr{3}{1}{3}{1.5}+\Teisssr{5}{2}{3}{1.5}$};
\draw[->] (0,2) -- (13,2);
\draw[->] (0,2) -- (0,6);
\draw[-, line width=0.5mm] (0,5) -- (5,3);
\draw[-, line width=0.5mm] (5,3) -- (8,2);
\end{scope}
\end{tikzpicture}
\end{center}
 \caption{Symbolic derivatives}
\label{fig:symbder}
    \end{figure}
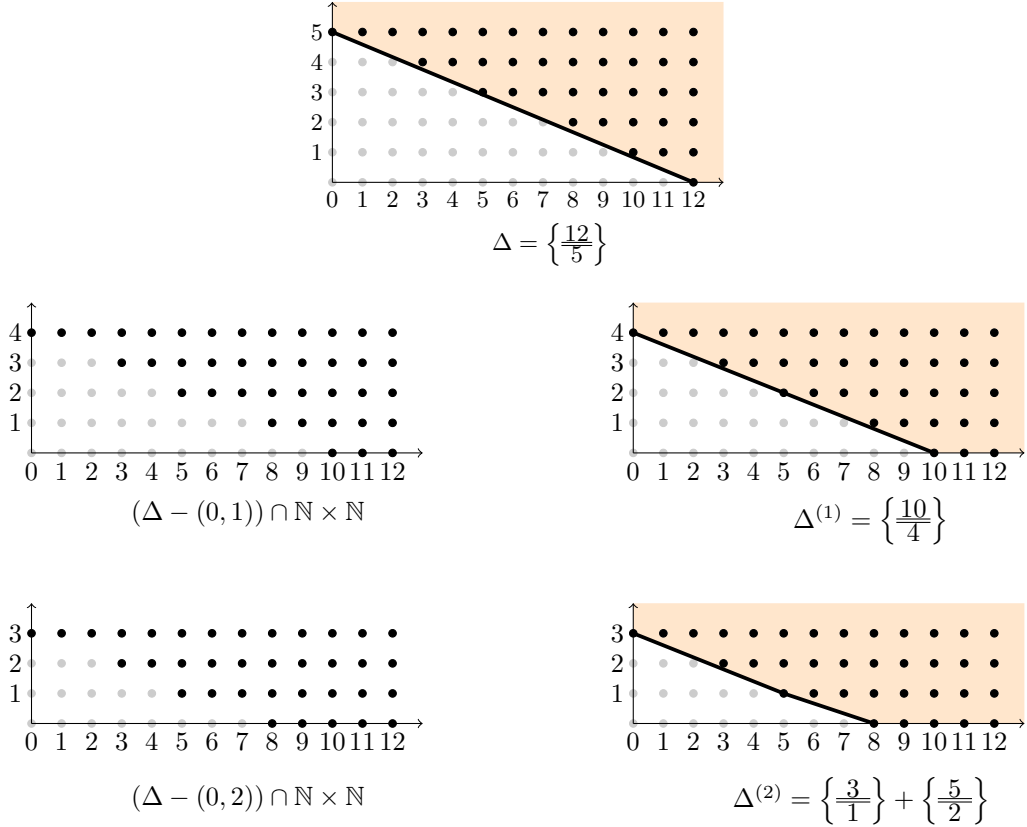

\begin{Example}
The symbolic first derivative of $\Delta=\Teisssr{12}{5}{3}{1.5}$ is $\Delta^{(1)}=\Teisssr{10}{4}{3}{1.5}$ and its symbolic second derivative is $\Delta^{(2)}=\Teisssr{3}{1}{3}{1.5}+\Teisssr{5}{2}{3}{1.5}$ (see Figure \ref{fig:symbder}).
\end{Example}

The main result of this paper is

\begin{teorema}\label{mainth} 
Let $f\in \mathbb C[[x,y]]$ be a generic element of $K(b_0,\ldots,b_h)$. 
Put $e_{i} =\gcd(b_{0}, \ldots , b_{i})$, $n_{i}=\frac{e_{i-1}}{e_{i}}$, $m_{i}=\frac{b_{i}}{e_{i}}$ 
and $\Delta_{i}=\Teisssr{m_{i}}{n_{i}}{4}{2}$ for $i\in \{1,\ldots,h\}$. Fix $1\leq  k < b_0$
and let $\{1,\dots,i_k\}=\{\,j\in\{1,\dots,h\}: e_{j-1}>k\,\}$.
Then $\frac{\partial^{k}f}{\partial y^k}$ admits the following factorization:
\[\frac{\partial^{k}f}{\partial y^k}=\Gamma^{(1)}\cdots \Gamma^{(i_{k})},\] 
where, for any $\ell\in \{1,\ldots,i_{k}\}$, the power series $\Gamma^{(\ell)}$ is not necessarily irreducible,  
and it verifies:
\begin{enumerate}
\item $\cont(f,v)=\frac{b_{\ell}}{b_{0}}$ for any irreducible factor $v$ of $\Gamma^{(\ell)}$.
\item Let $t$ be the natural number such that $0< t\leq n_\ell$ and  $t\equiv k$ mod $n_{\ell}$. If   $\sum_{j=1}^{r}\Teisssr{M_{j}}{N_{j}}{4}{2}$ is the long canonical representation of $\Delta_{\ell}^{(t)}$ and $m=\min\{e_{\ell}, k\}-\lceil\frac{k}{n_{\ell}}\rceil$  then $\Gamma^{(\ell)}$ can be written as a product of irreducible factors
\[\Gamma^{(\ell)}=\prod_{j=1}^{r}z^{(\ell)}_j\prod_{i=1}^{m}w^{(\ell)}_i \]
such that
\begin{enumerate}

\item[(2a)] for any power series $z^{(\ell)}_j$, $\,\cont(f_{\ell},z^{(\ell)}_j)=\frac{M_{j}}{n_{1}\cdots n_{\ell-1}N_{j}}$ and  

\[\Char(z^{(\ell)}_j)=
\left\{
\begin{array}{lr}
\left\{\frac{b_{1}}{b_{0}},\ldots, \frac{b_{\ell-1}}{b_{0}}\right\} & \hbox{\rm if } N_j=1\\
& \\
\left\{\frac{b_{1}}{b_{0}},\ldots, \frac{b_{\ell-1}}{b_{0}}, \frac{M_{j}}{n_{1}\cdots n_{\ell-1}N_{j}}\right\} & \hbox{\rm if } N_j>1.
\end{array}
\right . 
\]

\item[(2b)] for any power series $w^{(\ell)}_i$, $\,\Char(w^{(\ell)}_i)=\left\{\frac{b_{1}}{b_{0}},\ldots, \frac{b_{\ell}}{b_{0}}\right\}$ and   $\cont(f_{\ell},w^{(\ell)}_i)=\frac{b_{\ell}}{b_{0}}$.
\item[(2c)] $\cont(v_1,v_2)=\min\{\cont(f_l,v_1),\cont(f_l,v_2)\}$ for any two different irreducible factors $v_1,v_2$  of $\Gamma^{(\ell)}$.
\end{enumerate}
\end{enumerate}
\end{teorema}

Theorem \ref{mainth} will be proved in Section \ref{sec:proof}. It generalizes the results of  Casas-Alvero (see \cite{Casas}) 
and Hefez-Hernandes-Hern\'andez (see \cite{H-H-H}) to higher order polars.

\section{Symbolic derivatives of a Newton diagram}
In this section we  prove some properties of the symbolic derivatives of a Newton diagram. 

\begin{property} \label{trans} 
For any Newton diagram $\Delta$ and any nonnegative integers $k,l$ we have $(\Delta^{(k)})^{(l)}=\Delta^{(k+l)}$.
\end{property}
\begin{proof}
Note that 
$\supp(\Delta^{(k+l)})=(\supp(\Delta)-(0,k+l))\cap \mathbb N^2=(\supp(\Delta^{(k)})-(0,l))\cap \mathbb N^2=\supp(\Delta^{(k)})^{(l)}$.
\end{proof}

Let $\omega\in (\mathbb R_{> 0})^2$ and  $\Delta$ any Newton diagram. The {\em $\omega$-weigthed initial part} of $\Delta$ is 

\[
\ini_{\omega}(\Delta):=\{d\in \Delta\;:\;\langle d, \omega\rangle=\min\{\langle e ,\omega \rangle\;:\;e\in \Delta\}\},
\]

\noindent where $\langle \cdot  ,\cdot \rangle $ denotes the canonical scalar product in $\mathbb R^2$.

The Minkowski sum of Newton diagrams verifies the following  property (see \cite[Theorem 1.5, Chapter IV]{Ewald}).

\begin{property}\label{prop weight}
Let $\Delta_1, \Delta_2$ be two Newton diagrams and $\omega\in (\mathbb R_{> 0})^2$. Then
\[
\ini_{\omega}(\Delta_1+\Delta_2)=\ini_{\omega}(\Delta_1)+\ini_{\omega}(\Delta_2).
\]
\end{property}

\begin{nota}\label{vertices}
Let \begin{equation}\label{eq:cornc}\mathcal N=\sum_{i=1}^{r}\Teisssr{M_{i}}{N_{i}}{3}{1.5}\end{equation} be the canonical or the long canonical representation of $\mathcal N$. For any $0\leq j\leq r$ we put $A_j:=(a_j,b_j)$, where 
$a_j=\sum_{i=j+1}^rM_i$ and $b_j=\sum_{i=1}^jN_i$ (see Figure \ref{fig:points}).

 \begin{figure}[h!] 
  \begin{center}
\begin{tikzpicture}[x=0.7cm,y=0.7cm]
\tikzstyle{every node}=[font=\small]
\draw[->] (0,0) -- (5,0) node[right,below] {$\;$};
\draw[->] (0,0) -- (0,4) node[above,left] {$\;$};
\draw[thick] (4,1) node[right] {$A_{i-1}=(a_{i-1},b_{i-1})$};
\draw[thick] (1.2,3) node[above, right] {$A_i=(a_i,b_i)$};
\node[draw,circle,inner sep=1pt,fill] at (1,3) {};
\node[draw,circle,inner sep=1pt,fill] at (4,1) {};
\draw[dashed] (0,3) --(1,3);
\draw[dashed] (0,1) --(4,1);
\draw[dashed] (4,0) --(4,1);
\draw[dashed] (1,0) --(1,3);
\node[draw,circle,inner sep=1pt,fill] at (4,1) {};
\draw[ultra thick] (1,3) -- (4,1);
\draw[<->] (-0.5,3) to [bend right] (-0.5,1);
\draw[-] (-1.4,2.5) node[right,below] {$N_i$};
\draw[<->] (1,-0.5) to [bend right] (4,-0.5);
\draw[-] (2.7,-0.9) node[below] {$M_i$};
\end{tikzpicture}
\end{center}
\caption{ Points $A_{i-1}$ and $A_i$} 
\label{fig:points}
\end{figure}
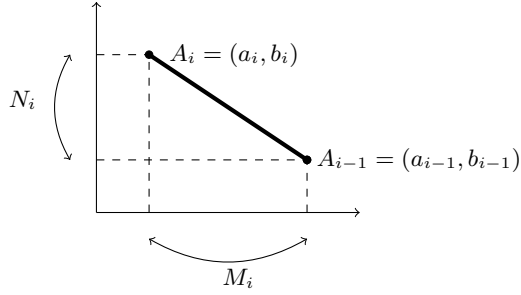

Observe that the set $T:=\{A_j\;:\; 0\leq j\leq r\}$ is a subset of the Newton polygon of $\mathcal N$ containing the vertices of $\mathcal N$,  with equality if and only if \eqref{eq:cornc} is the canonical representation of $\mathcal N$. In fact, if we consider $\omega_j:=(N_j,M_j)$ then $(M_i,0)\in \ini_{\omega_j}\left({\Teisssr{M_{i}}{N_{i}}{3}{1.5}}\right)$ 
for any $i>j$ 
and $(0,N_i)\in \ini_{\omega_j}\left(\Teisssr{M_{i}}{N_{i}}{3}{1.5}\right)$ 
for any $i \leq j$. 
By Property \ref{prop weight} we have $(a_j,b_j)\in \ini_{\omega_j}(\Delta)$, so $A_j\in \delta(\Delta)$.
\end{nota}

For any Newton polygon $\Delta$ let $\hbox{\rm trunc}(\Delta,k):=\mathcal{N}(\{(i,j)\in \Delta \cap \mathbb N^2\;:\; j\geq k\})$. It follows directly from definitions that
\begin{equation}\label{eq:twok}
\hbox{\rm trunc}(\Delta,k)=\Delta^{(k)}+(0,k).
\end{equation}

\begin{propo}\label{lr}
Let $\Delta=\sum_{i=1}^{r}\Teisssr{M_{i}}{N_{i}}{3}{1.5}$ be the long canonical representation of the convenient Newton diagram $\Delta$. 
Put $R=\sum_{i=1}^{s}\Teisssr{M_{i}}{N_{i}}{3}{1.5}$, 
$L=\sum_{i=s+1}^{r}\Teisssr{M_{i}}{N_{i}}{3}{1.5}$ 
and assume that $0\leq k\leq \sum_{i=1}^s N_i$. Then 
\begin{equation}\label{Dkk} 
\Delta^{(k)} = R^{(k)}+L. 
\end{equation}
\end{propo}
\begin{proof}
By Remark~\ref{vertices} the points $A_j=(a_j,b_j)=(\sum_{i=j+1}^rM_i,\sum_{i=1}^jN_i)$ 
for $0\leq j\leq r$ are  lattice points of  $\delta(\Delta)$. 

Let $R_1=\mathcal{N}(\{A_0,\dots,A_s\})$. Since the points 
$A_0,\dots,A_s$ are the lattice points 
of $\delta(R_1)$, we get by Remark~\ref{vertices} that 
$R_1=(a_s,0)+R$.    
Similarly for $L_1=\mathcal{N}(\{A_s,\dots,A_r\})$ we have $L_1=(0,b_s)+L$. 
Since $\Delta=L_1\cup R_1$, we get $\hbox{\rm trunc}(\Delta,k)=L_1\cup \hbox{\rm trunc}(R_1,k)$, 
and consequently $\hbox{\rm trunc}(\Delta,k)=L + \hbox{\rm trunc}(R,k)$.  Hence
equality \eqref{Dkk} follows.  
\end{proof}

\begin{corollary}
Let $\Delta=\sum_{i=1}^{r}\Teisssr{M_{i}}{N_{i}}{3}{1.5}$ be the long canonical representation 
of a convenient Newton diagram $\Delta$. Then
\[
\Delta^{(1)}=\Teisssr{M_{1}}{N_{1}}{3}{1.5}^{(1)}
+\sum_{i=2}^{r} \Teisssr{M_{i}}{N_{i}}{3}{1.5}.
\]
\end{corollary}

Recall the notion of continued fraction expansions of rational numbers.\\

Let $n,m\in \mathbb N$ with $0<n<m$. Denote by $[h_{0},h_{1},\ldots,h_{s}]$ the {\bf continued fraction expansion} of $\frac{m}{n}$, that is: 

\begin{equation}
\label{cf}
\frac{m}{n}=h_0+\cfrac{1}{h_1 +\cfrac{1}{h_2 +\cfrac{1}{\ddots + \cfrac{1}{h_s}}}}.
\end{equation}
Observe that the expansion \eqref{cf} is unique if we impose that $h_s>1$, that is, $s$ is the minimal possible value,  and this is the classical definition of continued fraction expansion. But if $h_s>1$ then $[h_{0},h_{1},\ldots,h_{s}]=
[h_{0},h_{1},\ldots,h_{s}-1,1]$. Hence, we can always assume, if necessary, that $s$ is  even.\\

Given the expansion \eqref{cf}, we put $p_{-1}=1$, $q_{-1}=0$, $p_{0}=h_{0}$, $q_{0}=1$ and consider the irreducible fractions
\[
\frac{p_{i}}{q_{i}}=[h_{0},h_{1},\ldots,h_{i}]=h_0+\cfrac{1}{h_1 +\cfrac{1}{h_2 +\cfrac{1}{\ddots + \cfrac{1}{h_i}}}}
\]

for $1\leq i\leq s$. The next properties are well-known (see for example \cite{Teissier2}).
\begin{properties}
\label{CF}
With the above notations we have:
\begin{enumerate}
\item $p_{i+1}=h_{i+1}p_{i}+p_{i-1}$ and $q_{i+1}=h_{i+1}q_{i}+q_{i-1}$, for $0\leq i\leq s-1$.
\item $p_{i}q_{i-1}-p_{i-1}q_{i}=(-1)^{i+1}$.
\item $\gcd(p_{i},q_{i})=1$.
\item $\frac{p_0}{q_0}<\frac{p_2}{q_2}<\cdots \leq \frac{m}{n}$.
\item $\frac{p_1}{q_1}>\frac{p_{3}}{q_{3}}>\cdots \geq \frac{m}{n}.$
\end{enumerate}
\end{properties}

Observe that $\frac{p_{s}}{q_{s}}=\frac{m}{n}$. If $m$ and $n$ are coprime then
$m=p_{s}$ and $n=q_{s}$. \\

\begin{propo}\label{Prop:derivative}
$\;$

 If $\Delta=\Teisssr{m}{n}{3}{1.5}$ with $n,m\in \mathbb N$ coprime then
\[
\Delta^{(1)}=\left\{\begin{array}{lr}\sum_{i=1}^{s/2}h_{2i}\Teisssr{p_{2i-1}}{q_{2i-1}}{6}{3}& \hbox{\rm if $s$ is even}\\
 \sum_{i=1}^{(s-1)/2}h_{2i}\Teisssr{p_{2i-1}}{q_{2i-1}}{6}{3}+\Teisssr{p_{s}-p_{s-1}}{q_{s}-q_{s-1}}{9}{4.5}& \hbox{\rm if $s$ is odd}.
 \end{array}
\right.
\]
In particular if $\Delta=\Teisssr{m}{1}{3}{1.5}$ then  $\Delta^{(1)}=\Teisssr{0}{0}{3}{1.5}$, that is the first quadrant.

\end{propo}
\begin{proof}
Suppose that $s$ is even. Consider the Newton diagram  
\[
N:=(0,1)+\sum_{i=1}^{s/2}h_{2i}\Teisssr{p_{2i-1}}{q_{2i-1}}{6}{3}=(0,1)+\sum_{i=1}^{s/2} \Teisssr{h_{2i}p_{2i-1}}{h_{2i}q_{2i-1}}{8}{4}.
\]

Since $p_1/q_1>p_3/q_3>\dots >p_{s-1}/q_{s-1}$, we get by Remark \ref{vertices} that the points 
\[
B_i:=\Bigl(\sum_{j=i+1}^{s/2}h_{2j}p_{2j-1}, \, 1 + \sum_{j=1}^{i}h_{2j}q_{2j-1}\Bigr)
\]
are the vertices of $\Delta(N)$ for 
$i=0, \dots, s/2$. 
By the first item of Properties \ref{CF} we get 
\[
B_i = \Bigl(\sum_{j=i+1}^{s/2}(p_{2j}-p_{2j-2}), \, 1 + \sum_{j=1}^{i}(q_{2j}-q_{2j-2})\Bigr)
    = (p_{s}-p_{2i}, q_{2i})
\]
for $i=0, \dots, s/2$.

We claim that
\begin{equation}\label{*}N=\cN(\{B_0, \dots, B_{s/2}\})=\hbox{\rm trunc}(\Delta,1).
\end{equation}

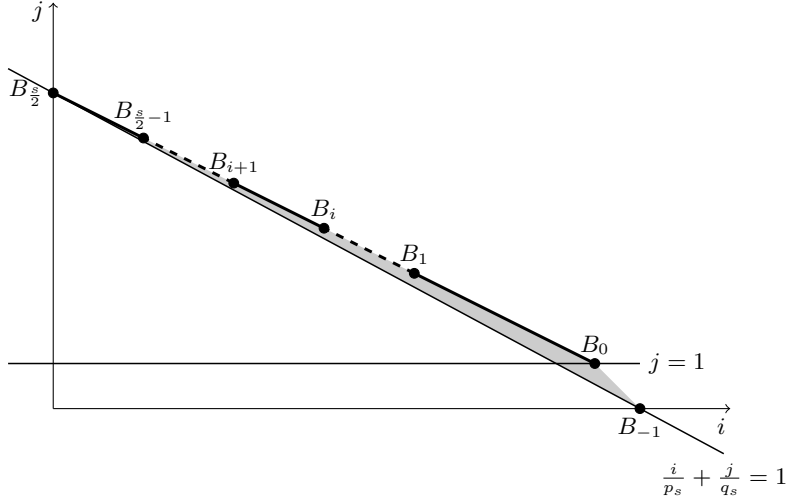
\begin{figure}[h!] 
\begin{center}
\begin{tikzpicture}[x=0.6cm,y=0.6cm] 
\tikzstyle{every node}=[font=\small]
\fill[fill=black!20!white] (0,7) --(12,1)--(13,0)--cycle;

\draw[->] (0,0) -- (15,0);
\node[left,below] at (14.8,0) {$i$};
\draw[->] (0,0) -- (0,9);
\node[below,left] at (0,8.8) {$j$};
\node[draw,circle,inner sep=1.3pt,fill, color=black] at (0,7){};
\node[draw,circle,inner sep=1.3pt,fill, color=black] at (2,6){};
\node[draw,circle,inner sep=1.3pt,fill, color=black] at (4,5){};
\node[draw,circle,inner sep=1.3pt,fill, color=black] at (6,4){};
\node[draw,circle,inner sep=1.3pt,fill, color=black] at (8,3){};
\node[draw,circle,inner sep=1.3pt,fill, color=black] at (12,1){};
\node[draw,circle,inner sep=1.3pt,fill, color=black] at (13,0){};
\draw[-, line width=0.2mm] (-1,98/13) -- (104/7,-1);
\draw[-, line width=0.2mm] (-1,1) -- (13,1);
\draw[-, line width=0.4mm] (12,1) -- (8,3);
\draw[-, dashed, line width=0.4mm] (8,3) -- (6,4);
\draw[-, line width=0.4mm] (6,4) -- (4,5);
\draw[-, dashed, line width=0.4mm] (4,5) -- (2,6);
\draw[-, line width=0.4mm] (2,6) -- (0,7);

\node [right] at (13,1){$j=1$};
\node [below] at (13,0){$B_{-1}$};
\node [above] at (12,1){$B_{0}$};
\node [above] at (8,3){$B_{1}$};
\node [above] at (6,4){$B_{i}$};
\node [above] at (4,5){$B_{i+1}$};
\node [above] at (2,6){$B_{\frac{s}{2}-1}$};
\node [left] at (0,7){$B_{\frac{s}{2}}$};
\node [below] at (104/7,-1){$\frac{i}{p_{s}}+\frac{j}{q_{s}}=1$};
\end{tikzpicture}
\end{center}
\caption{Points $B_i$}
\label{MMMM}
    \end{figure}

Consider the closed polygon $\mathcal{P}$ which vertices are $B_{-1}:=(p_s,0), B_0, \dots, B_{s/2}$ (see Figure \ref{MMMM}).

In order to proof equality \eqref{*} it is enough to show that there are no lattice points in the interior of $\mathcal{P}$. Denote by $B$ the number of lattice points belonging to the boundary of $\mathcal{P}$ and by $I$ the number of lattice points lying in the interior of $\mathcal{P}$. By the third item of Properties \ref{CF} we get $B=2+\sum_{i=0}^{s/2-1}h_{2i+2}$. By Pick's Formula \cite[Theorem 13.51]{Pick}, we have $2\mbox{Area}\,\mathcal{P}=2I+B-2$. On the other hand  if $\triangle_i$ denotes the triangle of vertices
$O,B_{i-1},B_i$ for 
$i=0, \dots, r$ then 
$2\mbox{Area}\,\mathcal{P}=\sum_{i=0}^{s/2} 2\mbox{Area}\,\triangle_i - p_sq_s$.
We have $2\mbox{Area}\,\triangle_0=p_s$ and 
$2\mbox{Area}\,\triangle_i=(p_s-p_{2i-2})q_{2i}-(p_s-p_{2i})q_{2i-2}=p_sq_{2i}-p_sq_{2i-2}+
p_{2i}q_{2i-2}-p_{2i-2}q_{2i}=p_sq_{2i}-p_sq_{2i-2}+(h_{2i}p_{2i-1}+p_{2i-2})q_{2i-2}-p_{2i-2}(h_{2i}q_{2i-1}+q_{2i-2})=
p_sq_{2i}-p_sq_{2i-2}+h_{2i}(p_{2i-1}q_{2i-2}-p_{2i-2}q_{2i-1})=
p_{s}q_{2i}-p_{s}q_{2i-2}+h_{2i}$ for $i=1,\dots,s/2$.  Hence 
\[
2\mbox{Area}\,\mathcal{P}= p_s+\sum_{i=1}^{s/2} (p_{s}q_{2i}-p_{s}q_{2i-2}+h_{2i})  - p_sq_s = 
\sum_{i=1}^{s/2} h_{2i}.
\]


Therefore $2\mbox{Area}\,\mathcal{P}=B-2$, and so $I=0$.

Suppose now that $s$ is odd. Observe that $[h_0,\ldots,h_s]$ can be represented by the even continued fraction $[h_0,\ldots,h_{s}-1,1]$. We have $[h_0,\ldots,h_{s}-1]=\frac{\tilde{p}_s}{\tilde{q}_s}$, where, by the first item of Properties \ref{CF}, we get $\tilde{p}_s=(h_s-1)p_{s-1}+p_{s-2}=p_s-p_{s-1}$ and $\tilde{q}_s=(h_s-1)q_{s-1}+q_{s-2}=q_s-q_{s-1}$. Hence the proof of $(1)$ for the odd case follows from the statement for the even case.
\end{proof}

\section{Technical tools}

The {\it extreme right edge} of a Newton polygon is its compact edge of greatest inclination.

\begin{lemma} \label{L:sliding}
Let $\lambda(x)=\sum_{i=1}^{N} a_ix^{i/n}$ be a finite Puiseux power series of characteristic 
$(n, b_1,\ldots, b_l)$ and let $v\in \C[[x,y]]$ be an irreducible power series 
such that $v(x,\lambda(x))\neq0$. 
Let 
\begin{equation}\label{cr}\sum_{i=1}^{r}\Teisssr{M_{i}}{N_{i}}{4}{2},
\end{equation}
be the canonical representation of the Newton diagram of 
$\hat v(x,y):=v(x^{n},y+\lambda(x^{n}))$. Then 

\begin{itemize}
\item[(i)] $M_i/N_i\leq N $ for all $i$ such that $1< i\leq r$,
\item[(ii)] $M_1/N_1= n\,\cont(\lambda,v)$,
\item[(iii)] if  $M_1/N_1>N$ and $\hat v$ is non-degenerate on the extreme right edge of its Newton polygon, 
then $M_1$ and $N_1$ are coprime and 
\[\Char(v)=
\left\{
\begin{array}{lr}
\left(\frac{b_{1}}{n},\ldots, \frac{b_{\ell}}{n}\right) & \hbox{\rm if } N_1=1\\
& \\
\left(\frac{b_{1}}{n},\ldots, \frac{b_{\ell}}{n}, \frac{M_{1}}{nN_{1}}\right) & \hbox{\rm if } N_1>1.
\end{array}
\right . 
\]
\end{itemize}
\end{lemma}

\begin{proof}
Let $\alpha_1,\dots,\alpha_m$ be the Newton-Puiseux roots of $v$. 
Then the set of Newton-Puiseux roots of $\hat v$ equals 
 $\{ \alpha_i(x^n)-\lambda(x^n): 1\leq i\leq m \} . $
Hence the set of inclinations of the edges of the Newton diagram 
of~$\hat v$ is equal to 
$ \{\, n\,\ord (\alpha_i(x)-\lambda(x)): 1\leq i\leq m \,\}$. 
In particular the biggest inclination $M_1/N_1$ of the Newton polygon of $\hat v(x,y)$
equals $n\,\cont(\lambda,v)$, 
which gives (ii). 

If $M_1/N_1\leq N$ then (i) clearly holds.  Hence in what follows, assume that $n\,\cont(\lambda,v)=M_1/N_1>N$. Then any Newton-Puiseux root $\alpha_i$
of $v$ that realizes the contact with $\lambda$ has the form 
$\alpha_i=\lambda+c_ix^{M_1/(nN_1)}+\cdots$ with some $c_i\neq0$.  
Thus for any $1\leq j\leq m$: 
either $\ord(\alpha_j-\lambda)=M_1/(n N_1)$ or 
$\ord(\alpha_j-\lambda)\leq N/n$. This proves (i).  

Assume that $\hat v$ is non-degenerate on the 
compact edge $S$ of \eqref{cr} of inclination $M_1/N_1$ and  suppose to the contrary that $\alpha_i$ has a characteristic exponent $\gamma$ bigger than $M_1/(n N_1)$. Then there exists $k\neq i$ such that $\gamma=\ord(\alpha_k-\alpha_i)$. This implies that 
$c_ix^{M_1/N_1}$ is the initial term of $\alpha_i(x^n)-\lambda(x^n)$ and 
$\alpha_k(x^n)-\lambda(x^n)$. 
In consequence, after (ND3), $\hat v$ is degenerate on the edge $S$ which is a  contradiction. Thus all characteristic 
exponens of  $\alpha_i$ are less than of equal to $M_1/(n N_1)$.

By \eqref{cr} there are $N_1$ Newton-Puiseux roots of $\hat v$ of order $\frac{M_1}{ N_1}$. Write $\frac{M_1}{n N_1}=\frac{m_{l+1}}{n\cdot n_{l+1}}$ with coprime 
$m_{l+1}$, $n_{l+1}$. After \eqref{eq:1} there are $n_{l+1}$ Newton-Puiseux roots $\alpha_j$
of $v$ such that $\ord(\alpha_j-\alpha_i) > \frac{b_{l}}{n}$. These Newton-Puiseux roots of $v$ yield the Newton-Puiseux roots 
$\alpha_j(x^n)-\lambda(x^n)$ of $\hat v$ of order $M_1/N_1$. 
Hence $N_1=n_{l+1}$. This proves (iii).
\end{proof}

\begin{corollary}\label{Coro}
Let $\lambda=\sum_{i=1}^N a_ix^{i/n}$ be a finite Puiseux series of characteristic 
$(n, b_1,\ldots, b_l)$ with minimal polynomial $g\in \C[[x]][y]$. 
Let $v\in \C[[x,y]]$ be a power series  coprime with $g$. 
Set $\hat v(x,y)=v(x^{n},y+\lambda(x^{n}))$. Let 
$$\sum_{i=1}^{s}\Teisssr{M_{i}}{N_{i}}{4}{2}$$ 
be the long canonical representation of 
${\mathcal N}(\hat v)$.
Assume that for some rational number $q\geq N$ the power series $\hat v$
is non-degenerate on all edges of inclination bigger than~$q$. 
Let $r$ be the number of elements of the set  
$\{\,i\in\{1,\dots,s\}: M_i/N_i>q\,\}$.
Then there exists a decomposition $v=\prod_{i=1}^a v_i$ into irreducible factors in $\C[[x,y]]$ such that:
\begin{itemize}
\item[(i)] $\cont(v_i,g)>q/n$ if and only if $1\leq i\leq r$,
\item[(ii)] for every $1\leq i\leq r$; $\cont(v_i,g)=\frac{M_{i}}{n\cdot N_{i}}$ 
 and 
 \[\Char(v_i)=
\left\{
\begin{array}{lr}
\left(\frac{b_{1}}{n},\ldots, \frac{b_{\ell}}{n}\right) & \hbox{\rm if } N_{i}=1\\
& \\
\left(\frac{b_{1}}{n},\ldots, \frac{b_{\ell}}{n}, \frac{M_{i}}{nN_{i}}\right) & \hbox{\rm if } N_i>1.
\end{array}
\right . 
\]
 
\item[(iii)] for every $1\leq i<j\leq r$;  $\cont(v_i,v_j)=\min\{\cont(v_i,g),\cont(v_j,g)\}$.
\end{itemize}
\end{corollary}

\begin{proof}
Let $v=\prod_{i=1}^a v_i$ be a decomposition of $v$ 
into irreducible factors in $\C[[x,y]]$ such that $\cont(g,v_i)\geq \cont(g,v_{i+1})$, for $1\leq i<a$.
Choose $r'\in \{1,\ldots,a\}$ such that $\cont(v_i,g)>q/n$ for any $1\leq i\leq r'$ and $\cont(v_i,g)\leq q/n$ for any $r'+1\leq i\leq a$.
Then by Lemma~\ref{L:sliding}, for $1\leq i\leq r'$, the Newton diagram of 
$\hat v_i:=v_i(x^{n},y+\lambda(x^{n}))$ has  one edge $L$ of inclination 
$n\,\cont(v_i,g)$ and all other edges have inclinations not greater than $q$.
Let $V_i=\hat v/\hat v_i$. 
Then $\mathcal{N}(\hat v)= \mathcal{N}(V_i)+\mathcal{N}(\hat v_i)$. 
In particular $\mathcal{N}(\hat v)$ has an edge $S$ of inclination 
$n\,\cont(v_i,g)$. Since $\hat v$ is non-degenerate on $S$ and 
the initial part of $\hat v_i$ with respect to $L$  divides the initial part of $\hat v$ with respect to $S$, we get, by $\hbox{\rm (ND1)}$ that 
$\hat v_i$ is non-degenerate on $L$. 
By~(ii) and~(iii) of Lemma~\ref{L:sliding}, the long canonical decomposition of ${\mathcal N}(\hat v_i)$  has only one  elementary Newton diagram (corresponding to $L$) of inclination greater than $q$. For $r'+1\leq i\leq a$, by~(i) and~(ii) of Lemma~\ref{L:sliding}, the inclinations of ${\mathcal N}(\hat v_i)$ are less than or equal to $q$.
\vspace{1cm}

From the identity 
\[ \mathcal{N}(\hat v) = \sum_{i=1}^a \mathcal{N}(\hat v_i) \]
we have $r=r'$ and
the extreme right compact edges of 
$\mathcal{N}(\hat v_i)$ for $1\leq i\leq r$  are in one-to-one correspondence with the set of elementary Newton diagrams 
$\left\{\Teisssr{M_{i}}{N_{i}}{4}{2}\right\}_{i=1}^r$ of the long canonical decomposition of $\mathcal N(\hat v)$.
Then, by (iii) of Lemma~\ref{L:sliding}, (i) and (ii) hold true. 

Suppose that there exists $1\leq i<j\leq r$ such that  the conclusion 
of~(iii) does not hold. 
This is possible only if $\cont(v_i,g)=\cont(v_j,g)<\cont(v_i,v_j)$.  
Let $\alpha_{i_0}$ be a Newton-Puiseux root of $v_i$ such that 
$\ord(\alpha_{i_0}-\lambda)=\cont(v_i,g)$ and let 
$\alpha_{j_0}$ be a Newton-Puiseux root of $v_j$
such that $\ord(\alpha_{j_0}-\alpha_{i_0})=\cont(v_i,v_j)$.  
Then the Puiseux series $\alpha_{i_0}(x^n)-\lambda(x^n)$, 
$\alpha_{j_0}(x^n)-\lambda(x^n)$ have the same initial term of 
order $n\,\cont(v_i,g)$. Hence, by \hbox{\rm (ND3)}, $\hat v$ is degenerate on the edge 
of inclination $n\,\cont(v_i,g)$. This contradiction gives (iii). 
\end{proof}

\begin{nota} \label{epi} For any positive  integers $r,s$   we have the epimorphism of groups $\mathbb U_r\ni\epsilon  \longrightarrow  \epsilon^s\in \mathbb U_{r/\gcd(r,s)}$. This becomes an isomorphism when $r, s$ are coprime.
\end{nota}

\medskip

\begin{properties} \label{propU}
Let $n\in \mathbb N$, $n>1$. Consider the strictly decreasing sequence $n=e_0>e_1>\cdots>e_h=1$ from page \pageref{page1}. Put $n_{i}=\frac{e_{i-1}}{e_{i}}$ for $1\leq i\leq h$. Then
\begin{enumerate}
    \item $\prod_{\varepsilon\in \mathbb U_{e_{l-1}}} (t-c\varepsilon^{b_l})  
    = (t^{n_l}-c^{n_l})^{e_l}$ for any $c\in \mathbb C$.
    \item $\prod_{\varepsilon\in \mathbb U_{e_{i-1}}\setminus \mathbb U_{e_i}} (1-\varepsilon^{b_i}) = n_i^{e_i} $.
  \item $ \sum_{\varepsilon\in \mathbb U_{n_l}} \varepsilon^{i} = 
\left\{\begin{array}{ll}
n_l,& \mbox{ if $i \equiv 0 \pmod{n_l}$} \\
0, &  \mbox{otherwise}.
\end{array}\right.
$ 
\end{enumerate}
\end{properties}
\begin{proof}

By Remark \ref{epi} $U_{e_{l-1}}\ni \varepsilon\to \varepsilon^{b_l}\in U_{n_l}$ is an epimorphism of groups, so
\[ \prod_{\varepsilon\in U_{e_{l-1}}} (t-c\varepsilon^{b_l}) =  \prod_{\tau\in U_{n_l}} (t-c\tau)^{e_l} = (t^{n_l}-c^{n_l})^{e_l}.\]
 In order to proof (2) consider
$h(x):=\prod_{\tau\in U_{n_i}\setminus \{1\}} (x-\tau)$. 
We have $(x-1)h(x)=x^{n_i}-1$, hence
$h(x)+(x-1)h'(x)=\frac{d}{dx} (x^{n_i}-1)=n_ix^{n_i-1}$. Substituting $x=1$ we get $h(1)=n_i$ which gives 
\[ \prod_{\varepsilon\in U_{e_{i-1}}\setminus U_{e_i}} (1-\varepsilon^{b_i}) = 
   \prod_{\tau\in U_{n_i}\setminus \{1\}} (1-\tau)^{e_i} = n_i^{e_i} . 
\] 

Statement (3) follows from Remark \ref{epi}.

\end{proof}

Let $f(x,y)=\sum_{ij}a_{ij}x^iy^j\in \mathbb C[[x,y]$ and $\omega=(\omega_1,\omega_2)\in \mathbb Q_{>0}^2$. The {\it $\omega$-weighted order of} $f$ is $\ord_{\omega}(f)=\min \{\omega_1 i+ \omega_2 j\,:\; a_{ij}\neq 0\}$
and the {\it $\omega$-weighted initial form of} $f$ is $\ini_{\omega}(f)=\sum_{ij}a_{ij}x^iy^j$, where the sum runs over $(i,j)$ such that $\omega_1 i+ \omega_2 j=\ord_{\omega}(f)$.

\begin{lemma}\label{ND}
Let $f\in K(n,b_1,\ldots,b_h)$ and 
$\alpha=\sum_{i\geq n} a_ix^{i/n}$ be a Newton-Puiseux root of $f$. Let 
$\lambda=\sum_{i=n}^{b_{l}-1} a_ix^{i/n}$
and $\hat f(x,y)=f(x^{n/e_{l-1}},y+\lambda(x^{n/e_{l-1}}))$.
Let $\Delta$ be the Newton diagram of $\hat f$. 
Then, for $k<e_{l-1}$ we have $\Delta^{(k)}=R^{(t)}+L$ where $R=\Teisssr{m_{l}}{n_{l}}{3}{1.5}$ and $t$ is the remainder 
of the division of $k$ by $n_l$.  The inclination of every compact edge of $L$ is smaller than or equal to $m_l/n_l$ and the inclination of every compact edge of $R^{(t)}$ is bigger than $m_l/n_l$. Moreover if $f$ is a generic member of $K(n,\ldots,b_h)$ then 
\begin{equation}\label{NDkd} \cN\left( \tfrac{\partial^k \hat f}{\partial y^k}\right) = \Delta^{(k)} 
\end{equation}
and $\frac{\partial^k \hat f}{\partial y^k}$ is non-degenerate on all edges of
its Newton diagram with inclinations  bigger than $m_l/n_l$.
\end{lemma}
\begin{proof} Observe that $\lambda(x^{n/e_{l-1}})\in \mathbb C[[x]]$, so $\hat f(x,y)=f(x^{n/e_{l-1}},y+\lambda(x^{n/e_{l-1}}))$ is a formal power series in $\mathbb C[[x,y]]$. The set of Newton-Puiseux roots of $\hat f(x,y)$ is $\Zer \hat f=\{\alpha_{\epsilon}(x^{n/e_{l-1}})-\lambda(x^{n/e_{l-1}})\;:\;\epsilon \in \mathbb U_n\}$. Hence $\{\ord (\gamma)\;:\; \gamma \in \Zer \hat f\}=\left\{ \frac{b_j}{e_{l-1}}\;:\; j=1,\ldots,l \right\}$. In particular the biggest inclination of the Newton diagram ${\mathcal N}(\hat f)$ equals $\frac{m_l}{n_l}$. Denote by $S$ the compact edge of ${\mathcal N}(\hat f)$ of this inclination. If $g$ is the minimal polynomial of $\lambda(x)$ then $g$ is a $l$-semiroot  of $f$, that is, $g\in \mathbb C[[x]][y]$ is monic, irreducible, its  $y$-degree equals $n/e_{l-1}$ and the intersection multiplicity of $f$ and $g$ is $\bar b_l:=b_l+\sum_{i=1}^{l-1}\left(\frac{e_{i-1}-e_i}{e_{l-1}}\right)b_i$ (see \cite[Theorem 3.9 (a)]{Zariski}). Hence the vertex of $S$ living on the horizontal axis is $(\overline b_l,0)$ since  $i_0( \hat f,y)=\ord(f(x^{n/e_{l-1}},\lambda(x^{n/e_{l-1}}))=i_0( f,g)$. On the other hand the length of the vertical projection of $L$ equals the cardinality of the set
\[
\{\alpha_{\epsilon}\in \Zer f\;:\;\ord (\alpha_{\epsilon}(x^{n/e_{l-1}})-\lambda(x^{n/e_{l-1}}))=\tfrac{m_l}{n_l}\}=
\{\alpha_{\epsilon}\in \Zer f\;:\;\ord (\alpha_{\epsilon}-\alpha)\geq \tfrac{b_l}{n}\}
\]

which is, after \eqref{eq:1}, equal to $e_{l-1}$. 

Let $k=qn_l+t$ be the Euclidean division of $k$ by $n_l$. Then $\Delta=q\Teisssr{m_{l}}{n_{l}}{3}{1.5}+\Teisssr{m_{l}}{n_{l}}{3}{1.5}+L$, for some Newton diagram $L$ with inclinations less than or equal  $\frac{m_l}{n_l}$. Consequently by Property \ref{trans} and Proposition \ref{lr} we get
$\Delta^{(k)}=\left(\Delta^{(qn_l)}\right)^{(t)}=\left(\Teisssr{m_{l}}{n_{l}}{3}{1.5}+L\right)^{(t)}=\Teisssr{m_{l}}{n_{l}}{3}{1.5}^{(t)}+L$.

\medskip

Now we are going to prove the second part of the lemma.\\

Suppose first that $f$ is a Weierstrass polynomial, that is $f$ is as in \eqref{NPf} with $u(x,y)=1$.
Then
\begin{equation}\label{NPfhat}
\hat f(x,y)=\prod_{\epsilon\in \mathbb U_n}(y-(\alpha_{\epsilon}(x^{n/e_{l-1}})-\lambda(x^{n/e_{l-1}}))).
\end{equation}
Fix  $q\in\{1,\dots,n_l-1\}$ and let
$z_q:=a_{b_l+qe_l}$ 
be a coefficient of $\alpha$ treated as indeterminate.  Expand $\hat f$ as a polynomial in $z_q$ 
\begin{equation}\label{fhatpoly}
    \hat f= \hat f_{q,0}+\hat f_{q,1}z_q+\cdots+\hat f_{q,n}z_q^n. 
\end{equation}
Consider $\omega:=(1,m_l/n_l)$.

{\it Claim 1.} The $\omega-$weighted order of $\hat f$ is $\ord_{\omega}(\hat f)=\bar b_l$ and the $\omega-$weighted initial form of $\hat f$ is
\begin{equation}\label{eq:init1}
 {\rm in}_{\omega}\hat f = a x^b (y^{n_l}-a_{b_l}^{n_l}x^{m_l})^{e_l}
 \end{equation}
for some nonzero complex number $a$ and a nonnegative integer $b$.  \\

Indeed, after \eqref{NPfhat} ${\rm in}_{\omega}\hat f=\prod_{\epsilon\in \mathbb U_n}{\rm in}_{\omega}A_{\epsilon}$, where $A_{\epsilon}=y-(\alpha_{\epsilon}(x^{n/e_{l-1}})-\lambda(x^{n/e_{l-1}}))$. Notice that
\[
{\rm in}_{\omega}A_{\epsilon}=\left\{
\begin{array}{ll}
(1-\epsilon^{b_j})a_{b_j}x^{b_j/e_{l-1}},& \hbox{\rm if } \epsilon\in \mathbb U_{e_{j-1}}\backslash \mathbb U_{e_{j}}\;\hbox{\rm for } 1\leq j\leq l-1\\
& \\
y-a_{b_j}\epsilon^{b_l}x^{b_l/e_{l-1}} & \hbox{\rm if } \epsilon\in \mathbb U_{e_{l-1}}.
\end{array}
\right.
\]
Hence 
\begin{eqnarray*}{\rm in}_{\omega}\hat f&=&\left(\prod_{j=1}^{l-1}\prod_{\epsilon \in \mathbb U_{e_{j-1}}\backslash \mathbb U_{e_{j}}}(1-\epsilon^{b_j})a_{b_j}x^{b_j/e_{l-1}}\right)\prod_{\epsilon \in \mathbb U_{e_{l-1}}}\left(y-a_{b_l}\epsilon^{b_l}x^{b_l/e_{l-1}}\right).
\end{eqnarray*}
By Properties \ref{propU} we get 
\begin{eqnarray*}{\rm in}_{\omega}\hat f
&=&\left(\prod_{j=1}^{l-1}n_j^{e_j}a_{b_j}^{e_{j-1}-e_j}x^{b_j(e_{j-1}-e_j)/e_{l-1}}\right)(y^{n_l}-a_{b_l}^{n_l}x^{b_l/e_l})^{e_l}.\\
\end{eqnarray*}
Notice that $b_l/e_l=m_l$ and the proof of  Claim 1 follows taking $a:=\prod_{j=1}^{l-1} n_j^{e_j}a_{b_j}^{e_{j-1}-e_j}$ and $b:=\ord_x \hat f(x,0)-b_l=\bar b_l-b_l\in \mathbb N$.\\

{\it Claim 2.} Let $q\in\{1,\dots,n_l-1\}$ and $\hat f_{q,0}$, $\hat f_{q,1}$  be as in \eqref{fhatpoly}. Then $\ord_{\omega}(\hat f_{q,1})=\bar b_l+\frac{q}{n_l}$ and
 \[{\rm in}_{\omega}(\hat f-\hat f_{q,0}) = {\rm in}_{\omega}\hat f_{q,1}z_q =
   -e_{l-1} a a_{b_l}^{s-1} x^{b+(m_ls+q)/n_l} y^{n_l-s} (y^{n_l}-a_{b_l}^{n_l}x^{m_l})^{e_l-1} z_q,
 \] where $s\in\{1,\dots,n_l\}$ 
is the solution of the congruence 
$m_ls+q\equiv 0 \pmod{n_l}$.

Indeed by Leibnitz rule 
\[\frac{d}{dz_q}\hat f=\hat f\sum_{\varepsilon\in \mathbb U_n} \frac{-\varepsilon^{b_l+qe_l} x^{(m_l+q)/n_l}}{A_{\varepsilon}}.\]

Hence by Remark \ref{epi} (for $r=e_{l-1}$ and $s=e_l$) we get
\begin{eqnarray} \label{eq:eq4}
{\rm in}_{\omega}\frac{d}{dz_q}\hat f&=&{\rm in}_{\omega}\hat f \cdot \left(\sum_{\varepsilon\in \mathbb U_{e_{l-1}}} \frac{-\varepsilon^{b_l+qe_l} x^{(m_l+q)/n_l}}{y-\varepsilon^{b_l}a_{b_l}x^{m_l/n_l}} \right) \nonumber\\
&= &
 {\rm in}_{\omega}\hat f \cdot \left(  \sum_{\varepsilon\in \mathbb U_{e_{l-1}}} \frac{-(\varepsilon^{e_l})^{m_l+q} x^{(m_l+q)/n_l}}{y-(\varepsilon^{e_l})^{m_l}a_{b_l}x^{m_l/n_l}} 
 \right) \nonumber\\
 &=&
 - e_lx^{(m_l+q)/n_l}{\rm in}_{\omega}\hat f \cdot\sum_{\theta\in \mathbb U_{n_l}} \frac{\theta^{m_l+q} }{y-\theta^{m_l}a_{b_l}x^{m_l/n_l}}.
\end{eqnarray} 

Let $q'$ be a solution of the congruence $m_lq'\equiv q\pmod{n_l}$. Then 
\begin{equation}\label{eq:eq3} \sum_{\theta\in \mathbb U_{n_l}} \frac{\theta^{m_l+q}}{y-\theta^{m_l}a_{b_l}x^{m_l/n_l}} = 
   \sum_{\theta\in \mathbb U_{n_l}} \frac{(\theta^{m_l})^{1+q'}}{y-\theta^{m_l}a_{b_l}x^{m_l/n_l}} =
   \sum_{\varepsilon\in \mathbb U_{n_l}} \frac{\varepsilon^{1+q'}}{y-\varepsilon a_{b_l}x^{m_l/n_l}}, 
\end{equation} 
where the last equality follows from Remark \ref{epi} for $r=n_{l}$ and $s=m_l$.\\

Using the equality 
$$ \frac{y^{n_l}-a_{b_l}^{n_l}x^{m_l}}{y-\varepsilon a_{b_l}x^{m_l/n_l}}=
   \sum_{j=0}^{n_l-1} \varepsilon^j a_{b_l}^jx^{jm_l/n_l} y^{n_l-1-j} $$
   for any $\varepsilon\in \mathbb U_{n_l}$ we have 

\begin{eqnarray}
\label{eq:eq2}
(y^{n_l}-a_{b_l}^{n_l}x^{m_l}) \sum_{\varepsilon\in \mathbb U_{n_l}} \frac{\varepsilon^{1+q'}}{y-\varepsilon a_{b_l}x^{m_l/n_l}} &=& 
\sum_{j=0}^{n_l-1} \sum_{\varepsilon\in \mathbb U_{n_l}} \varepsilon^{1+q'+j} a_{b_l}^jx^{jm_l/n_l} y^{n_l-1-j} \nonumber \\
&=&n_l a_{b_l}^{j_0}x^{j_0m_l/n_l} y^{n_l-1-j_0} 
\end{eqnarray}
where the last equality follows from the third part of Properties \ref{propU} 
and $j_0\in\{0,\dots,n_l-1\}$ satisfies $1+q'+j_0 \equiv 0 \pmod{n_l}$, that is $j_0$ is the solution of the congruence $m_l(j+1)+q\equiv 0 \pmod{n_l}$.\\

From \eqref{eq:eq4}, \eqref{eq:eq3} and \eqref{eq:eq2} it  follows  
\begin{equation}\label{eq:eq5} {\rm in}_{\omega}\tfrac{d}{dz_q}\hat f = \frac{{\rm in}_{\omega}\hat f}{(y^{n_l}-a_{b_l}^{n_l}x^{m_l})}(-1)e_{l-1} a_{b_l}^{s-1}x^{(m_ls+q)/n_l} y^{n_l-s} 
\end{equation}
with $s=j_0+1$.\\
We see that ${\rm in}_{\omega}\tfrac{d}{dz_q}\hat f$ does not depend on $z_q$. Thus, in view of the equality 
$\tfrac{d}{dz_q} \hat f = \hat f_{q,1} +2\hat f_{q,2} z_q+\cdots +n\hat f_{q,n} z_q^{n-1}$ we have 
${\rm in}_{\omega}\tfrac{d}{dz_q}\hat f={\rm in}_{\omega}\hat f_{q,1}$ and 
${\rm ord}_{\omega} (\hat f_{q,1})<{\rm ord}_{\omega} (\hat f_{q,j})$ for $j>1$. Consequently 
${\rm in}_{\omega}(\hat f-\hat f_{q,0})= 
{\rm in}_{\omega}(\hat f_{q,1}z_q + \cdots +\hat f_{q,n} z_q^{n}) = 
{\rm in}_{\omega}\hat f_{q,1}z_q$.
 Claim 2 follows from Claim 1 and \eqref{eq:eq5}.\\

 \noindent {\it Claim 3.} 
Let $q\in\{1,\dots,n_l-1\}$. Consider $u(x,y)\in \mathbb C[[x,y]]$, $u(0,0)=1$. Put $\hat u=u(x^{n/e_{l-1}},y+\lambda(x^{n/e_{l-1}})$. Then $\hat u \hat f$ is a polynomial in $z_q$ equal to $\hat u \hat f_{q,0}+\hat u \hat f_{q,1}z_q+\cdots+
\hat u \hat f_{q,n}z_q^n$, where  $\hat f_{q,i}$ is as in \eqref{fhatpoly}. Moreover
${\rm in}_{\omega}\hat u \hat f ={\rm in}_{\omega}\hat f$ and ${\rm in}_{\omega}(\hat u \hat f-\hat u\hat{f}_{q,0}) ={\rm in}_{\omega}(\hat f -\hat{f}_{q,0}$). \\

 \medskip

 Since $z_q=a_{b_l+qe_l}$ is not a coefficient of $\lambda(x)$ then $\hat u$ is independent of $z_q$.  The first part of the claim follows.  The second part also follows since ${\rm in}_{\omega}(\hat u)=1$ and the weighted initial part of a product is the product of the weighted initial parts of the factors.

 \medskip

 Consider now the truncation 
 $\hbox{\rm trunc}(\Delta,k)$
 and the lines $ L_q:i+\frac{m_l}{n_l}j=\bar b_l+\frac{q}{n_l}$ 
 where $q$ is a natural number verifying $0\leq q\leq n_l$.\\ 
 

  \noindent {\it Claim 4.} The lattice points on the compact edges of $\hbox{\rm trunc}(\Delta,k)$ with inclinations strictly bigger than $\frac{m_l}{n_l}$  belong to the lines $L_q$ with $0\leq q\leq n_l-1$.\\

  Indeed, consider $D:=\{(i,j)\in \mathbb R^2\;:\;\bar b_l\leq i+\frac{m_l}{n_l}j<\bar b_l+1\}\cap \{(i,j)\in \mathbb R^2\;:\;0\leq j\leq e_{l-1}\}$ (see Figure \ref{fig:D}). Observe that any lattice point $(i_0,j_0)$ in $D$ belongs to $\bigcup_{q=0}^{n_l-1}L_q$ since the rational number  $i_0+\frac{m_l}{n_l}j_0$ belonging to the interval $[\bar b_l,\bar b_l+1)$ has the form $\bar b_l+\frac{q}{n_l}$ for some $q\in\{0,\ldots,n_l-1\}$. Let $k<e_{l-1}$ and consider $d:=\min\{i\in \mathbb N\;:\;i+k\frac{m_l}{n_l}\geq \bar b_l\}$.

  Let ${\mathcal B}$ be the intersection of the compact edges of $\hbox{\rm trunc}(\Delta,k)$ and the strip $\mathbb R\times [k,e_{l-1}]$. Since $\hbox{\rm trunc}(\Delta,k)$ is contained in  $\Delta$ then ${\mathcal B}$ also. The set ${\mathcal B}$ is the graph of a piecewise linear, convex, decreasing function, contained in $L_0^+:=\{(i,j)\in \mathbb R^2\;:\;i+\frac{m_l}{n_l}j\geq \bar b_l\}$. The endpoints of ${\mathcal B}$ are $(b,e_{l-1})$, $(d,k)$. By convexity,  ${\mathcal B}$ is contained in $L_{n_l}^-:=\{(i,j)\in \mathbb R^2\;:\;i+\frac{m_l}{n_l}j< \bar b_l+1\}$
  so ${\mathcal B}\subseteq D$ and Claim 4 follows.\\
   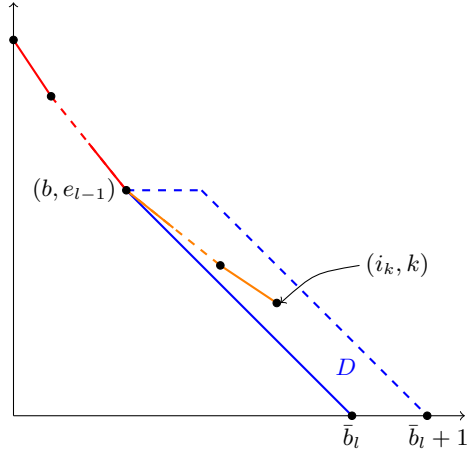
\begin{figure}[h!] 
  \begin{center}
\begin{tikzpicture}[x=0.5cm,y=0.5cm]
\tikzstyle{every node}=[font=\small]
\draw[->] (0,0) -- (12,0) node[right,below] {$\;$};
\draw[->] (0,0) -- (0,11) node[above,left] {$\;$};
\draw[thick] (9,0) node[below] {$\bar b_{l}$};
\draw[thick] (11.3,0) node[below] {$\bar b_{l}+1$};
\draw[thick] (3,6) node[left] {$(b,e_{l-1})$};
\draw[thick] (9.1,4) node[right]  {$(i_{k},k)$};
\draw[thick] (8.3,1.3) node[right,color=blue]  {$D$};

\draw[thick][-, color=red](0,10)--(1,8.5);
\draw[thick][dashed,color=red] (1,8.5) --(3,6);
\draw[thick][-,color=red] (2,7.25) --(3,6);
\draw[thick][-, color=blue](3,6)--(9,0);
\draw[thick][dashed,color=blue] (5,6) --(11,0);
\draw[thick][dashed,color=blue] (3,6) --(5,6);
\draw[thick][dashed,color=orange] (3,6) --(5.5,4);
\draw[thick][color=orange] (3,6) --(4.25,5);
\draw[thick][color=orange] (5.5,4)--(7,3);

\draw[->][color=black](9.2,4) .. controls (8,3.8) .. (7.1,3);
\node[draw,circle,inner sep=1pt,fill] at (0,10) {};
\node[draw,circle,inner sep=1pt,fill] at (1,8.5) {};
\node[draw,circle,inner sep=1pt,fill] at (3,6) {};
\node[draw,circle,inner sep=1pt,fill] at (9,0) {};
\node[draw,circle,inner sep=1pt,fill] at (11,0) {};
\node[draw,circle,inner sep=1pt,fill] at (5.5,4) {};
\node[draw,circle,inner sep=1pt,fill] at (7,3) {};
\end{tikzpicture}
\end{center}
\caption{The set D} 
\label{fig:D}
\end{figure}

Fix $q\in \{1,\ldots, n_l-1\}$. The lattice points of $L_q\cap D$ are the solutions of the linear diophantine equation $n_li+m_lj=n_l\bar b_l+q$ for $0\leq j\leq e_{l-1}$. Reducing this diophantine equation modulo $n_l$ we realize that there is no solution for $j=0$ or $j=e_{l-1}$. Hence the number of these lattice points is $e_l-1$  since $e_{l-1}=n_le_l$. Under the assumptions of Claim 3, the polynomial ${\rm in}_{\omega}(\hat u \hat f_{q,1})$ has $e_l-1$ monomials of $\omega$-weigthed order $\bar b_l+\frac{q}{n_l}$ and $y$-degree strictly less than $e_{l-1}$. Consequently these lattice points are in the support of ${\rm in}_{\omega}(\hat u \hat f_{q,1})$.\\

Let $f$ be a generic member of $K(n,b_1,\ldots,b_h)$. As the multiplication by nonzero constant does not affect the statement of the lemma, we may assume that $f=uf^*$ where $f^*\in \mathbb C[[x]][y]$ is a Weierstrass polynomial and $u(x,y)\in \mathbb C[[x,y]]$ with $u(0,0)=1$. \\

Observe that equality  \eqref{NDkd} is equivalent to equality $\cN\left( y^k\tfrac{\partial^k \hat f}{\partial y^k}\right) = \hbox{\rm trunc}(\Delta,k)$. Moreover it follows from (ND2) that $\frac{\partial^k \hat f}{\partial y^k}$ is non-degenerate on the edge $S$ of
its Newton diagram  if and only if $y^k\frac{\partial^k \hat f}{\partial y^k}$ is also non-degenerate on the edge $S+(0,k)$ of its Newton diagram.

Let $\{(i_r,j_r)\}_{r=0}^s$ be the set of lattice points belonging to the compact edges of ${\mathcal B}$ with inclinations strictly bigger than $\frac{m_l}{n_l}$, ordered by the first coordinate, that is $i_r<i_{r+1}$ for any $r\in \{0,\ldots, s-1\}$. Note that $(i_s,j_s)=(d,k)$. We have $\bar b_l=\ord_{\omega}(x^{i_0}y^{j_0})<\ord_{\omega}(x^{i_1}y^{j_1})<\cdots <\ord_{\omega}(x^{i_{s}}y^{j_{s}})<\bar b_l+1$. For any $r\in\{1,\ldots, s\}$ there exists $q_r\in \{1,\ldots, n_l-1\}$ such that $\ord_{\omega}(x^{i_r}y^{j_r})=\bar b_l+\frac{q_r}{n_l}$.

Put $y^k\tfrac{\partial^k \hat f}{\partial y^k}=\sum c_{ij}x^iy^j$.\\

By Claims 3 and 1,  for any $r\in\{1,\ldots, s\}$, we have
 \begin{equation}\label{eqND}
     c_{i_rj_r}=W_{r}(z_1,\ldots,z_{q_r-1})+\gamma_rz_{q_r},
 \end{equation}
 where $\gamma_r\in \mathbb C\backslash \{0\}$, $W_r\in \mathbb C[z_1,\ldots,z_{q_r-1}]$ and 
 $c_{i_0j_0}$ is a nonzero constant polynomial in $\mathbb C[z_1,\ldots,z_{q_s}]$. The map  
 \[
 \begin{array}{rll}\Phi:\mathbb C^s& \longrightarrow & \mathbb C^s\\
 (z_{q_1},\ldots,z_{q_s})& \longrightarrow& \Phi(z_{q_1},\ldots,z_{q_s})=(c_{i_1j_1}(z_{q_1},\ldots,z_{q_s}), \ldots, c_{i_sj_s}(z_{q_1},\ldots,z_{q_s}))
 \end{array}
 \]
 is surjective after the triangular form of its components given by \eqref{eqND}. The equality
 $\cN\left( y^k\tfrac{\partial^k \hat f}{\partial y^k}\right) = \hbox{\rm trunc}(\Delta,k)$ is equivalent to the non-vanishing of all coefficients  $c_{ij}$ where $(i,j)\in \{(i_r,j_r)\}_{r=1}^s$ is a vertex of $\mathcal B$.\\
 
Assume for a moment that the equality $\cN\left( y^k\tfrac{\partial^k \hat f}{\partial y^k}\right) = \hbox{\rm trunc}(\Delta,k)$ holds. Let $R$ be a compact edge of  $\hbox{\rm trunc}(\Delta,k)$ of inclination bigger than $\frac{m_l}{n_l}$. Denote by $\alpha_R$ the maximum natural number $i$ such that $y^{i}$ divides the initial form  $g_R$ of  $g:=y^k\tfrac{\partial^k \hat f}{\partial y^k}$ with respect to $R$. 
 The non-degeneracy of $y^k\tfrac{\partial^k \hat f}{\partial y^k}$ on the compact edge $R$   is equivalent to the non-vanishing of the discriminant of the polynomial $y^{-\alpha_R}g_R(1,y)$.  Denote by $H_R$ this discriminant. Since the coefficients of $y^{-\alpha_R}g_R(1,y)$ are in the set $\{c_{i_{q_\ell}j_{q_\ell}}\}_{\ell=0}^s$ then $H_R\in \mathbb C[c_{i_{q_1}j_{q_1}},\ldots, c_{i_{q_s}j_{q_s}}]\backslash\{0\}$.

 Consider 
 \[{\mathcal A}_1:=\{c_{ij}\;:\;(i,j)\in \{(i_r,j_r)\}_{r=1}^s \;\hbox{\rm is a vertex of $\mathcal B$}\}\]
 
 and 
 
 \[
 {\mathcal A}_2:=\left\{H_R\;:\;R \;\hbox{\rm is a compact edge of }\hbox{\rm trunc}(\Delta,k)\; \hbox{\rm  of inclination bigger than }\tfrac{m_l}{n_l}\right \}.
 \]

The complement of the solutions of the polynomial defined as the product of all elements of  ${\mathcal A}_1 \cup {\mathcal A}_2$ is a non-empty open Zariski set in the target of $\Phi$ and its preimage by $\Phi$ is a non-empty open Zariski set in the source of $\Phi$. Hence there is a non-empty open Zariski  set in the space of coefficients of the Puiseux root $\alpha(x)$ of $f\in K(n,b_1,\ldots,b_g)$  such that

\[ \cN\left( \tfrac{\partial^k \hat f}{\partial y^k}\right) = \hbox{\rm trunc}(\Delta,k)
\]
and $\frac{\partial^k \hat f}{\partial y^k}$ is non-degenerate on all edges of
its Newton diagram which inclinations are bigger than $m_l/n_l$. 
This last non-empty open Zariski is the complement of the solutions of a polynomial depending on a finite number of coefficients of $\alpha$, let us say  $a_{s_1},\ldots, a_{s_{\ell}}$; and we denote this polynomial by $G(a_{s_1},\ldots, a_{s_{\ell}})$. Consider now the polynomial $\overline G=\prod_{\epsilon \in \mathbb U_n}G(\epsilon^{s_1}a_{s_1},\ldots, \epsilon^{s_{\ell}}a_{s_{\ell}})$. By \cite[Theorem 3]{G-G}, there exists a finite set of coefficients of $f$, let us say $a_{u_1v_1},\ldots, a_{u_Iv_I}$ and a polynomial $W\in \mathbb C[T_1, \ldots, T_I]$ such that $W(a_{u_1v_1},\ldots, a_{u_Iv_I})=0$ if and only if  $\overline G(a_{s_1},\ldots, a_{s_{\ell}})=0$. We conclude that 
if $f$ is a generic element in $K(n,b_1,\ldots, b_h)$, that is  $W(a_{u_1v_1},\ldots, a_{u_Iv_I})\neq 0$, 
then $G(a_{s_1},\ldots, a_{s_{\ell}})\neq 0$ and the lemma follows.
\end{proof}

\section{Proof of the main theorem}

\label{sec:proof}

In this section we will prove Theorem \ref{mainth}.
Let $f$ be a generic member of $K(b_0,\ldots,b_h)$. Remember that $e_i=\gcd(b_0,\ldots,b_i)$, for $0\leq i\leq h$ and $n_i=\frac{e_{i-1}}{e_i}$, $m_i=\frac{b_i}{e_i}$, $\Delta_i=\Teisssr{m_{i}}{n_{i}}{3}{1.5}$ for $1\leq i\leq h$. 
Fix $1\leq  k < b_0$ and 
let $\ell \in \{1,\ldots, h\}$ be such that $e_{\ell-1}>k$. Let $\alpha$ be any Newton-Puiseux root of $f$.

Denote by  $\lambda_{\ell}$ the sum of all terms of $\alpha$ of degree strictly less than $\frac{b_{\ell}}{b_0}$ and denote by $f_{\ell}(y)\in \mathbb C[[x]][y]$ the minimal polynomial of $\lambda_{\ell}$. The degree of $f_{\ell}(y)$ equals $n_1\cdots n_{\ell-1}$. Let $\frac{\partial^k f}{\partial y^k}=g_1\cdots g_r$ be the factorization into irreducible factors of the $k$th derivative of $f$. Put $\Gamma_{\ell}:=\prod_jg_j$ where the product runs over the factors $g_j$ such that $\cont(g_j,f)=\frac{b_{\ell}}{b_0}$.  By \cite[Theorem 6.2 ]{Forum} $\frac{\partial^{k}f}{\partial y^k}=\Gamma^{(1)}\cdots \Gamma^{(i_{k})}$ which proves item (1) of the theorem.

We can write $\Gamma^{(\ell)}=\Gamma^{(\ell)}_1\Gamma^{(\ell)}_2$ verifying $\cont(g,f_{\ell})> \frac{b_{\ell}}{b_0}$ for any irreducible factor $g$ of $\Gamma^{(\ell)}_{1}$ and $\cont(g,f_{\ell})=\frac{b_{\ell}}{b_0}$ for any irreducible factor $g$ of $\Gamma^{(\ell)}_{2}$.

Remark that the factors $\Gamma^{(\ell)}_{1}$ and $\Gamma^{(\ell)}_{2}$ coincide with the factors given in~\cite[Theorem 6.2 ]{Forum}.

After \cite[Theorem 6.2 (v), (ii)]{Forum} $\Gamma^{(\ell)}_{2}=\prod_{i=1}^{m}w^{(\ell)}_i$ where $m=\min\{e_{\ell}, k\}-\lceil\frac{k}{n_{\ell}}\rceil$ and the set of characteristic exponents of its irreducible factors $w^{(\ell)}_i$
 is $\left\{\frac{b_{1}}{b_{0}},\ldots, \frac{b_{\ell}}{b_{0}}\right\}$.

Since $\frac{b_{\ell}}{b_{0}}$ is not in the support of $\lambda_\ell$ we get $\cont(f_\ell, w_i^{(\ell)})=b_\ell/b_0$ for $1\leq i\leq m$, and statement $(2b)$ follows.\\

On the other hand we get $\widehat{\frac{\partial^k f}{\partial y^k}}=\frac{\partial^k \hat{f}}{\partial y^k}$, so by Lemma \ref{ND}   $\cN\left( \tfrac{\partial^k \hat f}{\partial y^k}\right) =\Teisssr{m_{\ell}}{n_{\ell}}{3}{1.5}^{(t)}+L$, where the inclinations of the compact edges of $L$ are less than or equal to $\frac{m_\ell}{n_\ell}$. Moreover $\widehat{\frac{\partial^k f}{\partial y^k}}$ is non-degenerate on all edges of
its Newton diagram which inclinations are bigger than $m_l/n_l$. 

\medskip
Now applying Corollary \ref{Coro} to $\lambda=\lambda_{\ell}$, which characteristic  is $\left(\frac{b_{0}}{e_{\ell-1}},\frac{b_{1}}{e_{\ell-1}},\ldots, \frac{b_{\ell-1}}{e_{\ell-1}}\right)$, $g=f_l$, $v=\frac{\partial^k f}{\partial y^k}$ and $q=\frac{m_\ell}{n_\ell}$, we get that $\Gamma_1^{(\ell)}$ can be written as $\prod_{j=1}^{r}z^{(\ell)}_j$ with $z^{(\ell)}_j$ irreducible verifying statements $(2a)$ and $(2c)$ of the theorem. 

\medskip In order to prove statement $(2c)$ in full generality it is enough to show that
$\cont(w_i^{(\ell)},w_j^{(\ell)})=\frac{b_\ell}{b_0}$ for $1\leq i<j\leq m$. Suppose that $\cont(w_i^{(\ell)},w_j^{(\ell)})>\frac{b_\ell}{b_0}$ for some $i,j\in \{1,\ldots,m\}$, $i\neq j$. Then there is a
 nonzero complex number $u$ and a 
Newton-Puiseux root $\gamma_d$   of $w_d^{(\ell)}$  such that $\gamma_d=\lambda_\ell+ux^{b_\ell/b_0}+\cdots$, for $d=i,j$. We claim that $u$ is not a root of the univariate polynomial $y^{n_\ell}-a_{b_\ell}^{n_\ell}$.
Indeed suppose that $u=\tau^{n_\ell}a_{b_\ell}$ for some $n_\ell$-th root of unity $\tau$. Let $\varepsilon$ an $e_{\ell-1}$-th root of the unity such that $\tau=\varepsilon^{b_\ell}$. Then the Newton-Puiseux root $\alpha_{\varepsilon}$ of $f$ has the form $\alpha_{\varepsilon}=\lambda_\ell+\varepsilon^{b_\ell}a_{b_\ell}x^{b_\ell/b_0}+\cdots=\lambda_\ell+ux^{b_\ell/b_0}+\cdots$, hence $\ord(\alpha_\varepsilon-w_d^{(\ell)})>b_l/b_0$ which is a contradiction since $\cont(f,w_d^{(\ell)})=b_l/b_0$ and we finished the proof of the claim.

Observe that $\tilde \gamma_d:=\gamma_d(x^{n/e_{\ell-1}})-\lambda_\ell(x^{n/e_{\ell-1}})=ux^{m_\ell/n_l}+\cdots$ 
are Newton-Puiseux roots of $\frac{\partial^k \hat f}{\partial y^k}$, for $d=i,j$. 

Let $F(y):=\ini_\omega\hat f(x,y)\vert_{x=1}$ (see \eqref{eq:init1}). 
Hence we get 
$
\tfrac{d^kF}{dy^k}=\ini_\omega\left(\tfrac{\partial^k \hat f(x,y)}{\partial y^k}\right)_{\vert_{x=1}}.
$
Given that $(y-ux^{m_\ell/n_\ell})^2$ is a factor of  $\ini_\omega\left(\tfrac{\partial^k \hat f(x,y)}{\partial y^k}\right)$ then $(y-u)^2$ is a factor of $ \tfrac{d^kF}{dy^k}$
which is a contradiction since  $\tfrac{d^kF}{dy^k}$  has no multiple complex roots except $0$ and the roots of $y^{n_\ell}-a_{b_\ell}^{n_\ell}$(see \cite[Corollary 5.4]{Forum}).The proof of Theorem \ref{mainth} is finished.

\medskip

\begin{Example}\label{Ex1}
Consider a generic element $f$ of $K(12,16,31)$.
Then 
\[ \frac{\partial f}{\partial y}=\Gamma^{(1)}\Gamma^{(2)} \]
where
\[
\cont(f,v)=\left\{
\begin{array}{ll}
\frac{4}{3} & \text{for any irreducible factor $v$ of $\Gamma^{(1)}$}\\
 & \\
\frac{31}{12} & \text{for any irreducible factor $v$ of $\Gamma^{(2)}$.}\\
\end{array}
\right .
\]
We have $(n_1,m_1)=(3,4)$ and $(n_2,m_2)=(4,31)$. 
The first symbolic derivatives of Newton diagrams $\Delta_1=\Teisssr{4}{3}{3}{1.5}$, $\Delta_2=\Teisssr{31}{4}{3}{1.5}$ are
$\Delta_1^{(1)}=\Teisssr{3}{2}{3}{1.5}$, $\Delta_2^{(1)}=3\Teisssr{8}{1}{3}{1.5}$. 
Hence $\Gamma^{(1)}=z_1^{(1)}$ and $\Gamma^{(2)}=\prod_{j=1}^3 z_j^{(2)}$ where

\begin{itemize}
\item $\cont(f_1,z_1^{(1)})=\frac{3}{2}$ and $\Char(z^{(1)}_1)=\left\{\frac{3}{2} \right \}$,
\item $\cont(f_2,z_j^{(2)})=\frac{8}{3}$ and  $\Char(z^{(2)}_j)=\left\{\frac{4}{3} \right \}$, for $j\in\{1,2,3\}$.
\end{itemize}

For the second polar we have $\frac{\partial^2 f}{\partial y^2}=\Gamma^{(1)}\Gamma^{(2)}$ where as before

\[
\cont(f,v)=\left\{
\begin{array}{ll}
\frac{4}{3} & \text{for any irreducible factor $v$ of $\Gamma^{(1)}$}\\
 & \\
\frac{31}{12} & \text{for any irreducible factor $v$ of $\Gamma^{(2)}$.}\\
\end{array}
\right .
\]

Since 
$\Delta_1^{(2)}=\Teisssr{2}{1}{3}{1.5}$, $\Delta_2^{(2)}=2\Teisssr{8}{1}{3}{1.5}$, we get in this case that 
$\Gamma^{(1)}=z_1^{(1)}w_1^{(1)}$ and $\Gamma^{(2)}=z_1^{(2)} z_2^{(2)}$ where

\begin{itemize}
\item $\cont(f_1,z_1^{(1)})=\frac{2}{1}$ and $\Char(z^{(1)}_1)=\emptyset$, that is, $z^{(1)}_1$ is smooth,

\item $\cont(f_1,w_1^{(1)})=\frac{4}{3}$ and  $\Char(w^{(1)}_1)=\left\{\frac{4}{3} \right \}$,

\item $\cont(f_2,z_j^{(2)})=\frac{8}{3}$ and  $\Char(z^{(2)}_j)=\left\{\frac{4}{3} \right \}$, for $j\in\{1,2\}$.
\end{itemize}

\medskip

Consider now $g(x,y)\in K(12,16,31)$ which admits $\alpha(x)=x^{4/3}+x^2+x^{31/12}$ as a Newton-Puiseux root. 
Applying a symbolic computation program {\bf Maxima} we get
$g(x,y)=y^{12}-12x^{2}y^{11}+66x^{4}y^{10}+h(x,y)$, where $\deg_yh(x,y)=9$. Hence $\frac{\partial^{10}g}{\partial y^{10}}=6\cdot 11!(y-x^2)^2$. However, after Theorem \ref{mainth}, for a generic element $f\in K(12,16,31)$ we get,
$\frac{\partial^{10}f}{\partial y^{10}}=\Gamma^{(1)}=z_1^{(1)}$, with $\Char(z^{(1)}_1)=\left\{\frac{3}{2} \right \}$, $\cont(f,z_1^{(1)})=\frac{4}{3}$ and $\cont(f_1,z_1^{(1)})=\frac{3}{2}$. 
We conclude that $g$ is not a generic element of $K(12, 16, 31)$ in the sense of Theorem \ref{mainth},

\end{Example}

\begin{Example}\label{Ex2}
Consider a generic element $f$ of $K(10,14,15)$. 
We have $\Delta_1=\Teisssr{m_1}{n_1}{3}{1.5}=\Teisssr{7}{5}{3}{1.5}$ and $\Delta_2=\Teisssr{m_2}{n_2}{3}{1.5}=\Teisssr{15}{2}{3}{1.5}$.
By Proposition~\ref{Prop:derivative} the first symbolic derivatives of these Newton diagrams are $\Delta_1^{(1)}=2\Teisssr{3}{2}{3}{1.5}$ and 
$\Delta_2^{(1)}=\Teisssr{8}{1}{3}{1.5}$. 

We get \[\frac{\partial f}{\partial y}=\Gamma^{(1)}\Gamma^{(2)},\]
where
\[
\cont(f,v)=\left\{
\begin{array}{ll}
\frac{7}{5} & \text{for any irreducible factor $v$ of $\Gamma^{(1)}$}\\
 & \\
\frac{3}{2} & \text{for any irreducible factor $v$ of $\Gamma^{(2)}$.}\\
\end{array}
\right .
\]
Moreover  $\Gamma^{(1)}=z_1^{(1)}z_2^{(1)}$
and $\Gamma^{(2)}= z_1^{(2)}$ where

\begin{itemize}
\item $\cont(f_1,z_j^{(1)})=\frac{3}{2}$ and $\Char(z^{(1)}_j)=\left\{\frac{3}{2} \right \}$, for $j\in\{1,2\}$;
\item $\cont(f_2,z_1^{(2)})=\frac{8}{5}$ and  $\Char(z^{(2)}_1)=\left\{\frac{7}{5} \right \}$.
\end{itemize}

For the second polar we have $\frac{\partial^2 f}{\partial y^2}=\Gamma^{(1)}$ where $\cont(f,v)=\frac{7}{5}$ for any irreducible factor $v$ of $\Gamma^{(1)}$. 

In this case $\Delta_1^{(2)}=\Teisssr{2}{1}{3}{1.5}+\Teisssr{3}{2}{3}{1.5}$. Hence $\Gamma^{(1)}=z_1^{(1)}z_2^{(1)}w_1^{(1)}$
where 

\begin{itemize}
\item $\cont(f_1,z_1^{(1)})=\frac{2}{1}$ and  $z^{(1)}_1$ is smooth,
\item $\cont(f_1,z_2^{(1)})=\frac{3}{2}$ and $\Char(z^{(1)}_2)=\{\frac{3}{2}\}$, 
\item $\cont(f_1,w_1^{(1)})=\frac{7}{5}$ and  $\Char(w^{(1)}_1)=\left\{\frac{7}{5} \right \}$. 
\end{itemize}
\end{Example}

\begin{Remark}

In Figures \ref{fig:EW1} and \ref{fig:EW2} we illustrate Examples \ref{Ex1} and~\ref{Ex2} using Eggers-Wall trees.
Recall that the Eggers-Wall tree $\Theta(h)$ of a reduced power series $h(x,y)$ is a rooted tree with leaves corresponding to irreducible factors of $h$. 
For any two irreducible factors $h_1,h_2$ of $h$ the last common vertex of the paths from the root of $\Theta(h)$ to $h_1$ and from the root to $h_2$ 
is labelled by the contact $\cont(h_1,h_2)$.  The Eggers-Wall tree $\Theta(h)$ equipped with some additional information (weights of edges)  
characterizes the equisingularity class of $h(x,y)$ (see \cite{Wall} and \cite{Handbook}). 

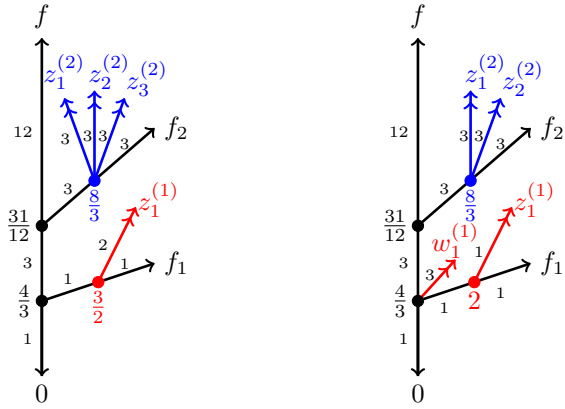
\begin{figure}
    \begin{center}
\begin{tikzpicture}[scale=1]
\begin{scope}[shift={(0,0)}]

     \draw [-, color=black,  line width=1pt](0,0) -- (0,2);
     \draw [->, color=black,  line width=1pt](0,2) -- (1.5,3.3);
      \draw [->, color=black,  line width=1pt](0,2) -- (0,4.5);
     \draw [->, color=black,  line width=1pt](0,1) -- (0,0);
     \draw [->, color=black,  line width=1pt](0,1) -- (1.5,1.5); 
     \draw [->_>, color=red,  line width=1pt](0.75,1.25) -- (1.25,2.25);
\draw [->_>, color=blue,  line width=1pt](0.7,2.6) -- (0.7,3.8);
\draw [->_>, color=blue,  line width=1pt](0.7,2.6) -- (0.3,3.7);
\draw [->_>, color=blue,  line width=1pt](0.7,2.6) -- (1.1,3.7);

       \node [below, color=black] at (0,0) {$0$};
       \node [above, color=black] at (0,4.5) {$f$};
           \node [right, color=black] at (1.5,3.3) {$f_2$};
              \node [right, color=black] at (1.5,1.5) {$f_{1}$};
              
            \node[draw,circle, inner sep=1.5pt,color=black, fill=black] at (0,1){};
      \node[draw,circle, inner sep=1.5pt,color=black, fill=black] at (0,2){};
           \node [left, color=black] at (0,1) {$\frac{4}{3}$}; 
            \node [left, color=black] at (0,2) {$\frac{31}{12}$}; 
            \node[draw,circle, inner sep=1.5pt,color=red, fill=red] at (0.75,1.25){};
\node[draw,circle, inner sep=1.5pt,color=blue, fill=blue] at (0.7,2.6){};

            \node [below, color=red] at (0.75,1.25) {$\frac{3}{2}$};
            \node [left] at (0,0.5) {\tiny{$\textcolor{black}{1}$}};
            \node [left] at (0,1.5) {\tiny{ $\textcolor{black}{3}$}};
            \node [left] at (0,3.25) {\tiny {$12$}};
            \node [above,right] at (1.15,2.35) { $\textcolor{red}{z^{(1)}_1}$};

            \node [above] at (0.35,1.1) { \tiny{${1}$}};
            \node [above] at (1.1,1.3) { \tiny{${1}$}};
            \node [left] at (1,1.75) { \tiny{${2}$}};
            \node [below, color=blue] at (0.7,2.6) {$\frac{8}{3}$};
            \node [above] at (0.3,3.7) { $\textcolor{blue}{z^{(2)}_1}$};
            \node [above] at (0.9,3.7) { $\textcolor{blue}{z^{(2)}_2}$};
            \node [above] at (1.4,3.6) { $\textcolor{blue}{z^{(2)}_3}$};
            \node [above] at (0.35,2.3) { \tiny{${3}$}};
            \node [above] at (1.1,2.9) { \tiny{${3}$}};
            \node [left] at (0.5,3.15) { \tiny{${3}$}};
            \node [left] at (0.8,3.2) { \tiny{${3}$}};
            \node [left] at (1,3.2) { \tiny{${3}$}};
                   
      \end{scope}
      
      \begin{scope}[shift={(5,0)}]

     \draw [-, color=black,  line width=1pt](0,0) -- (0,2);
     \draw [->, color=black,  line width=1pt](0,2) -- (1.5,3.3);
      \draw [->, color=black,  line width=1pt](0,2) -- (0,4.5);
     \draw [->, color=black,  line width=1pt](0,1) -- (0,0);
     \draw [->, color=black,  line width=1pt](0,1) -- (1.5,1.5); 
     \draw [->_>, color=red,  line width=1pt](0.75,1.25) -- (1.25,2.25);
     \draw [->_>, color=red,  line width=1pt](0,1) -- (0.5,1.55);
\draw [->_>, color=blue,  line width=1pt](0.7,2.6) -- (0.7,3.8);
\draw [->_>, color=blue,  line width=1pt](0.7,2.6) -- (1.1,3.7);

       \node [below, color=black] at (0,0) {$0$};
       \node [above, color=black] at (0,4.5) {$f$};
           \node [right, color=black] at (1.5,3.3) {$f_2$};
              \node [right, color=black] at (1.5,1.5) {$f_{1}$};
              
            \node[draw,circle, inner sep=1.5pt,color=black, fill=black] at (0,1){};
      \node[draw,circle, inner sep=1.5pt,color=black, fill=black] at (0,2){};
           \node [left, color=black] at (0,1) {$\frac{4}{3}$}; 
            \node [left, color=black] at (0,2) {$\frac{31}{12}$}; 
            \node[draw,circle, inner sep=1.5pt,color=red, fill=red] at (0.75,1.25){};
\node[draw,circle, inner sep=1.5pt,color=blue, fill=blue] at (0.7,2.6){};

            \node [below, color=red] at (0.75,1.25) {$2$};
            \node [left] at (0,0.5) {\tiny{$\textcolor{black}{1}$}};
            \node [left] at (0,1.5) {\tiny{ $\textcolor{black}{3}$}};
            \node [left] at (0,3.25) {\tiny {$12$}};
            \node [above,right] at (1.15,2.35) { $\textcolor{red}{z^{(1)}_1}$};

\node [above] at (0.5,1.45) { $\textcolor{red}{w^{(1)}_1}$};
\node [left] at (0.35,1.35) {\tiny {$3$}};
            
            \node [below] at (0.35,1.1) { \tiny{${1}$}};
            \node [below] at (1.1,1.3) { \tiny{${1}$}};
            \node [left] at (1,1.65) { \tiny{${1}$}};
            \node [below, color=blue] at (0.7,2.6) {$\frac{8}{3}$};
            \node [above] at (0.9,3.7) { $\textcolor{blue}{z^{(2)}_1}$};
            \node [above] at (1.4,3.6) { $\textcolor{blue}{z^{(2)}_2}$};
            \node [above] at (0.35,2.3) { \tiny{${3}$}};
            \node [above] at (1.1,2.9) { \tiny{${3}$}};
            \node [left] at (0.8,3.2) { \tiny{${3}$}};
            \node [left] at (1,3.2) { \tiny{${3}$}}; 
      \end{scope}
  \end{tikzpicture}
\end{center}
 \caption{Eggers-Wall trees of Example \ref{Ex1}: on the left $\Theta(ff_1f_2\frac{\partial f}{\partial y})$ and on the right $\Theta(ff_1f_2\frac{\partial^2 f}{\partial y^2})$}
\label{fig:EW1}
   \end{figure}

   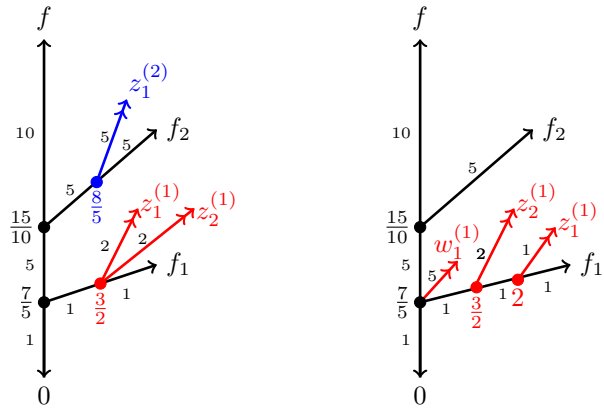
\begin{figure}
    \begin{center}
\begin{tikzpicture}[scale=1]
\begin{scope}[shift={(0,0)}]

     \draw [-, color=black,  line width=1pt](0,0) -- (0,2);
     \draw [->, color=black,  line width=1pt](0,2) -- (1.5,3.3);
      \draw [->, color=black,  line width=1pt](0,2) -- (0,4.5);
     \draw [->, color=black,  line width=1pt](0,1) -- (0,0);
     \draw [->, color=black,  line width=1pt](0,1) -- (1.5,1.5); 
     \draw [->_>, color=red,  line width=1pt](0.75,1.25) -- (1.25,2.25);
     \draw [->_>, color=red,  line width=1pt](0.75,1.25) -- (2,2.25);

\draw [->_>, color=blue,  line width=1pt](0.7,2.6) -- (1.1,3.7);

       \node [below, color=black] at (0,0) {$0$};
       \node [above, color=black] at (0,4.5) {$f$};
           \node [right, color=black] at (1.5,3.3) {$f_2$};
              \node [right, color=black] at (1.5,1.5) {$f_{1}$};
              
            \node[draw,circle, inner sep=1.5pt,color=black, fill=black] at (0,1){};
      \node[draw,circle, inner sep=1.5pt,color=black, fill=black] at (0,2){};
           \node [left, color=black] at (0,1) {$\frac{7}{5}$}; 
            \node [left, color=black] at (0,2) {$\frac{15}{10}$}; 
            \node[draw,circle, inner sep=1.5pt,color=red, fill=red] at (0.75,1.25){};
\node[draw,circle, inner sep=1.5pt,color=blue, fill=blue] at (0.7,2.6){};

            \node [below, color=red] at (0.75,1.25) {$\frac{3}{2}$};
            \node [left] at (0,0.5) {\tiny{$\textcolor{black}{1}$}};
            \node [left] at (0,1.5) {\tiny{ $\textcolor{black}{5}$}};
            \node [left] at (0,3.25) {\tiny {$10$}};
            \node [above,right] at (1.15,2.35) { $\textcolor{red}{z^{(1)}_1}$};
            \node [right] at (1.9,2.25) { $\textcolor{red}{z^{(1)}_2}$};

            \node [below] at (0.35,1.1) { \tiny{${1}$}};
            \node [below] at (1.1,1.3) { \tiny{${1}$}};
            \node [left] at (1,1.75) { \tiny{${2}$}};
            \node [left] at (1.5,1.85) { \tiny{${2}$}};
            \node [below, color=blue] at (0.7,2.6) {$\frac{8}{5}$};
            \node [above] at (1.4,3.6) { $\textcolor{blue}{z^{(2)}_1}$};
            \node [above] at (0.35,2.3) { \tiny{${5}$}};
            \node [above] at (1.1,2.9) { \tiny{${5}$}};
            \node [left] at (1,3.2) { \tiny{${5}$}};
                   
      \end{scope}
      
      \begin{scope}[shift={(5,0)}]

     \draw [-, color=black,  line width=1pt](0,0) -- (0,2);
     \draw [->, color=black,  line width=1pt](0,2) -- (1.5,3.3);
      \draw [->, color=black,  line width=1pt](0,2) -- (0,4.5);
     \draw [->, color=black,  line width=1pt](0,1) -- (0,0);
     \draw [->, color=black,  line width=1pt](0,1) -- (2,1.5); 
     \draw [->_>, color=red,  line width=1pt](0.75,1.25) -- (1.25,2.25);
     \draw [->_>, color=red,  line width=1pt](0,1) -- (0.5,1.55);
     \draw [->_>, color=red,  line width=1pt](1.3,1.3) -- (1.8,2);

       \node [below, color=black] at (0,0) {$0$};
       \node [above, color=black] at (0,4.5) {$f$};
           \node [right, color=black] at (1.5,3.3) {$f_2$};
              \node [right, color=black] at (2,1.5) {$f_{1}$};
              
            \node[draw,circle, inner sep=1.5pt,color=black, fill=black] at (0,1){};
      \node[draw,circle, inner sep=1.5pt,color=black, fill=black] at (0,2){};
           \node [left, color=black] at (0,1) {$\frac{7}{5}$}; 
            \node [left, color=black] at (0,2) {$\frac{15}{10}$}; 
            \node[draw,circle, inner sep=1.5pt,color=red, fill=red] at (0.75,1.2){};
            \node[draw,circle, inner sep=1.5pt,color=red, fill=red] at (1.3,1.3){};

            \node [below, color=red] at (1.3,1.3) {$2$};
            \node [below, color=red] at (0.75,1.2) {$\frac{3}{2}$};
            \node [left] at (0,0.5) {\tiny{$\textcolor{black}{1}$}};
            \node [left] at (0,1.5) {\tiny{ $\textcolor{black}{5}$}};
            \node [left] at (0,3.25) {\tiny {$10$}};
            \node [above,right] at (1.15,2.35) { $\textcolor{red}{z^{(1)}_2}$};
            \node [above, right] at (1.7,2.1) { $\textcolor{red}{z^{(1)}_1}$};

\node [above] at (0.5,1.45) { $\textcolor{red}{w^{(1)}_1}$};
\node [left] at (0.35,1.35) {\tiny {$5$}};
            
            \node [below] at (0.35,1.1) { \tiny{${1}$}};
            \node [below] at (1.1,1.3) { \tiny{${1}$}};
            \node [below] at (1.7,1.4) { \tiny{${1}$}};
            \node [left] at (1,1.65) { \tiny{${2}$}};
            \node [left] at (1.6,1.7) { \tiny{${1}$}};
            \node [left] at (1,1.65) { \tiny{${2}$}};
            \node [above, color=black] at (0.7,2.6) {\tiny{$5$}};
                  
      \end{scope}
  \end{tikzpicture}
\end{center}
 \caption{Eggers-Wall trees of Example \ref{Ex2}: on the left $\Theta(ff_1f_2\frac{\partial f}{\partial y})$ and on the right $\Theta(ff_1f_2\frac{\partial^2 f}{\partial y^2})$}
\label{fig:EW2}
   \end{figure}

\end{Remark}

\end{document}